\documentclass[12pt,twoside]{article}

\advance\textheight by 0.795in
\advance\textwidth by 1.19in \advance\oddsidemargin by -.595in
\advance\evensidemargin by -.595in
\usepackage{mathrsfs,upgreek,url}
\usepackage[all,cmtip]{xy}
\usepackage{amsfonts,amssymb,latexsym,amsthm}
\usepackage{epsfig,pstricks}
\usepackage{amsmath}
\usepackage{hyperref}
\usepackage{epsfig,pstricks}
\usepackage{amscd, enumerate}
\usepackage{verbatim}

\numberwithin{equation}{section}

\newtheorem{st}{statement}[section]
\newtheorem{lem}[st]{Lemma}

\newtheorem*{rep.thm}{\it{Theorem}}
\newtheorem{prop}[st]{Proposition}
\newtheorem{thm}[st]{Theorem}

\theoremstyle{definition}
\newtheorem{rep.lem}[st]{Lemma}
\newtheorem{para}[st]{}

\newtheorem*{con}{Notational Convention}
\newtheorem*{defn}{Definition}
\newtheorem*{i.thm}{Theorem}
\newtheorem*{notn}{Notation}

\newtheorem*{pf}{Proof}

\newtheorem*{key}{\it{Keywords}}

\newtheorem{rem}[st]{Remark}
\newtheorem*{Rem}{{Remark}}
\def\Bb{\mathbb}

\def\frk{\mathfrak}
\def\iff{if and only if\  }
\def\ocirc{\overset{\circ}}
\def\lb{\lbrace}
\def\rb{\rbrace}

\def\eps{\varepsilon}
\def\ad{{\text{ad}}}
\def\aV{Vogan }

\title{ Vogan Diagrams of Twisted\\ Affine Kac-Moody Lie 
Algebras}
\date{}
\author{Tanusree Pal\vspace{0.02cm}\\{\small Harish-Chandra Research Institute,
Chhatnag Road, Jhunsi,Allahabad -211019, India}\vspace{0.005cm}\\
{\small E-mail address :tanusree@hri.res.in}}
\begin{document}
\maketitle

\begin{abstract}
A Vogan diagram is a Dynkin diagram of a Kac-Moody Lie algebra of finite or 
affine type overlayed with additional structures. This paper develops the 
theory of \aV diagrams for ``almost compact'' real forms of indecomposable 
twisted affine Kac-Moody Lie algebras and shows that equivalence classes of 
Vogan diagrams correspond to isomorphism classes of almost compact real forms 
of twisted affine Kac-Moody Lie algebras as given by H. Ben Messaoud and 
G. Rousseau in the paper ``Classification des formes r\'{e}elles presque 
compactes des alg\`ebres de Kac Moody affines, J. Algebra 267". 
\end{abstract}

\begin{key}Almost compact real forms; Vogan diagram; Twisted affine Kac-Moody 
algebra. 
\end{key}
\maketitle

\section{Introduction}
The classification of finite dimensional real simple Lie algebras has been
a classical problem. In 1914, \'Elie Cartan classified the simple Lie algebras
over the reals for the first time. A number of subsequent simplifications of 
the proof followed, until in 1996, using the theory of Vogan diagrams,  
A. W. Knapp derived a quick proof of \'Elie Cartan's classification in 
\cite{KnappQ}. 

The Kac-Moody Lie algebras are an infinite-dimensional generalization of the 
semisimple Lie algebras via the Cartan matrix and generators. The real forms of
complex affine Lie algebras are of two kinds, ``almost split'' and 
``almost compact.'' V. Back et. al. classified the almost split real forms of 
the affine Kac-Moody Lie algebra  in \cite{1type} and  H. Ben Messaoud and 
G. Rousseau gave a classification of the almost compact real forms in 
\cite{M-R}. Work towards developing the theory of Vogan diagrams for the real 
forms of non-twisted affine Kac-Moody Lie algebras was done by P. Batra in 
\cite{B1, B2}. In the present paper the theory of Vogan diagrams for twisted
affine Kac-Moody Lie algebras is developed.

As introduced in \cite{Knapp}, a Vogan diagram is a Dynkin diagram of a Lie 
algebra with a diagram involution, such that the vertices fixed by the 
involution are either painted or unpainted depending on whether they are 
noncompact or compact. An important result in the theory of Vogan diagrams for
real simple Lie algebras states that any Vogan diagram can be transformed, 
by changing the ordering of its base, into a diagram which has at most one 
noncompact imaginary root and that root occurs at most twice in the largest 
root of that simple Lie algebra.  Since in the case of affine Kac-Moody 
algebras, changing the order does not give a Vogan diagram with at most one 
shaded root, therefore a notion of equivalence of Vogan diagrams for 
non-twisted affine Kac-Moody Lie algebras was introduced in \cite{B1}. In the 
present paper we modify the definition of the Vogan diagrams for the twisted 
affine Kac-Moody Lie algebras. In addition to the structural information 
already superimposed, a \aV diagram now contains numerical labels on the 
vertices of the underlying Dynkin diagram as given in Figure 1. The 
classification of the almost compact real forms of affine Kac-Moody Lie 
algebras as given in \cite{M-R}, prompts the definition of suitable equivalence
relations among the \aV diagrams for twisted affine Kac-Moody Lie algebras. 
With respect to this equivalence relation we prove the following result. 

\begin{i.thm} Let $\mathfrak{g}$ be a twisted affine Kac-Moody Lie algebra. 
Then
\begin{enumerate}
\item  Two almost compact real forms of $\frak{g}$ having equivalent 
Vogan diagrams are isomorphic.
\item Every abstract Vogan diagram for $\frak{g}$, represents an almost compact
real form of $\mathfrak{g}$. 
\end{enumerate}\end{i.thm}
The analogues of these results for the non-twisted Kac-Moody Lie algebras were 
proved by P. Batra in \cite[Theorem 5.2]{B1} and \cite[Theorem 5.2]{B2} 
respectively. Owing to the difference in the structural realizations of the 
non-twisted and twisted affine Kac-Moody Lie algebras, the methods used in 
\cite{B1, B2} prove insufficient to yeild the main theorems for the twisted 
affine Kac-Moody Lie algebras. This difficulty is resolved by using the notion of 
``adapted realization'' of an affine Kac-Moody Lie algebra as introduced in 
\cite{M-R}. 

It is a fairly easy matter to work out representatives of the equivalence 
classes of Vogan diagrams for the twisted affine Kac-Moody Lie algebras. These 
have been listed in Figures 2 and 3. A match in the count of the non-equivalent
Vogan diagrams and the count of non-isomorphic almost compact real forms as
given in \cite{M-R}, seems to suggest the existence of a bijective 
correspondence between the equivalence classes of Vogan diagrams and the 
isomorphism classes of almost compact real forms of twisted affine Kac-Moody 
Lie algebras.

This paper is organized as follows: Section 2 reviews known facts about 
(complex) indecomposable twisted affine Kac-Moody Lie algebras $\frk{g}$. In 
Section 3, the automorphisms and real forms of $\frk{g}$ are discussed, certain
results from \cite{M-R}, which allow simple proofs of the main theorems are 
recalled and some properties of the Cartan subalgebras of the almost compact 
real forms of $\frk{g}$ are studied. In Section 4, the \aV diagrams are 
introduced, their equivalence relations defined and the main theorems are 
stated and proved. In Figure 2 and 3, the non-equivalent \aV diagrams for
the twisted affine Kac-Moody Lie algebras are given.

\con The complexification $\frk{m}\otimes\Bb{C}$ of a real Lie algebra 
$\frk{m}$ will be denoted by $\frk{m}_{\Bb{C}}$. Given a finite order 
automorphism $\phi$ of a Lie algebra $\frk{L}$, we shall denote by 
$\frk{L}^{\phi}$, the fixed point subalgebra 
$\lbrace x\in \frk{L}\mid \phi(x)=x \rbrace$ of $\frk{L}$. By abuse of notation
we shall denote $\Bb{Z}\cap[m,n]$ by $[m,n]$ for all $m,n \in \Bb{Z}$. For all
integers $n$, $\varepsilon_n$ will denote the $n^{th}$ root of unity.

\section{Kac-Moody Lie algebras}
\para Let $\frk{g}$ be an affine Kac-Moody Lie algebra over the complex field 
$\Bb{C}$. There exists a generalized Cartan matrix $A=(a_{i,j})_{i,j\in [0,l]}$
such that $\frk{g}=\frk{g}(A)$ is generated by the Cartan subalgebra $\frk{h}$
and the elements $e_i,f_i$ for $i\in [0,l]$(cf. \cite[Chapter 1]{Kac}). We have
a decomposition $\frk{g}= \frk{h}\oplus(\bigoplus_{\alpha\in\triangle} 
\frk{g}_{\alpha})$, where $\triangle\subset{\frk{h}^*\smallsetminus \lb0\rb}$ 
denotes the root system of ($\frk{g},\frk{h}$). Let $\pi=\lb\alpha_i \mid i\in
[0,l] \rb$ be the standard base of $\triangle$; $\triangle_+ = \triangle\cap
(\oplus_{i\in[0,l]}\Bb{N}\alpha_i)$ the set of positive roots and
$\triangle_-=-\triangle_+$ the set of negative roots of $\frk{g}$. The coroots 
$(\alpha_i^{\vee})_{i\in[0,l]} \subset \frk{h}$ are such that 
$a_{i,j} = \alpha_j(\alpha_i^{\vee})$ for $i,j \in [0,l]$. Let $W$ denote the 
Weyl group of $\frk{g}$. $\alpha \in \triangle$ is said to be a real root, if
$\alpha$ is $W$-conjugate to a root in $\pi$ and the set of real roots is 
denoted by $\triangle^{re}$. The elements of $\triangle^{im} = 
\triangle\smallsetminus\triangle^{re}$ are called the imaginary roots of 
$\frk{g}$.

\para{\it{Realization of a Kac-Moody Lie algebra}} :
Let $\dot{\frk{g}}$ be a finite dimensional simple Lie algebra over $\Bb{C}$, 
$\mu$ a $k$-order automorphism of $\frk{\dot{g}}$ for $k < \infty$, $\eps_k = 
e^{\frac{2i\pi}{k}}$, a primitive $k^{th}$ root of unity and
$(.,.)$ a nondegenerate, invariant, symmetric bilinear form on $\dot{\frk{g}}$.
For $j\in\Bb{Z}_k$, let $\mathfrak{\dot{g}}_j$ = $\lbrace X \in \frk{\dot{g}} 
\mid \mu(X) = \varepsilon^{j}_k X \rbrace$; then $\frk{\dot{g}} =
\overset{k-1}{\underset{j=0}{\oplus}}\frk{\dot{g}}_j$. By \cite{Kac}, 
$(\frk{\dot{g}})^{\mu} = \ocirc{\frk{g}}=\frk{\dot{g}}_0$ is a simple finite 
dimensional Lie algebra over $\Bb{C}$. If $\dot{\frk{h}}$ denotes the Cartan 
subalgebra of $\frk{\dot{g}}$, then $\ocirc{\frk{h}} = 
\frk{\dot{h}}\cap\frk{\dot{g}}_0$ is the Cartan subalgebra of 
$\frk{\dot{g}}_0$. We denote by $\frk{g}$ the infinite dimensional Lie algebra:
$${\frk{l}}(\frk{\dot{g}},\mu,\eps_k) = \big(\bigoplus_{j\in\Bb{Z}} 
\frk{\dot{g}}_{(j\mod k)}\otimes t^j\big)\oplus\Bb{C}c\oplus \Bb{C}d.$$
The Lie algebra structure on $\frk{g}$ is such that $c$ is the canonical 
central element and
$$\lbrack  x\otimes t^m +\lambda d, y\otimes t^n + \lambda_{1} d \rbrack = 
([x,y] \otimes t^{m+n}+\lambda n y\otimes t^{n} -\lambda_{1}m x\otimes t^{m}) +
m\delta_{m,-n}(x,y) c, $$ where  $x,y \in 
\frk{\dot{g}}$, $\lambda,\lambda_{1}\in\Bb{C}$. The element $d$ acts 
diagonally on $\frk{g}$ with integer eigenvalues and induces $\Bb{Z}$-gradation
on ${\frk{l}}(\frk{\dot{g}},\mu,\eps_k)$. The Lie algebra 
$\frk{l}(\frk{\dot{g}},Id,1)$ with $\mu=Id$ denotes a non-twisted affine
Kac-Moody Lie algebra and $ {\frk{l}}(\frk{\dot{g}},\mu,\eps_k)$ for $\mu\neq 
Id$ and $k=2$ or 3 denotes a twisted affine Kac-Moody Lie algebra. 
Clearly $\frk{g}''= {\frk{l}}(\frk{\dot{g}},\mu,\eps_k)''= 
\underset{j\in\Bb{Z}}{\oplus} \frk{\dot{g}}_{(j\mod k)}\otimes t^j$ is the 
fixed point set of the automorphism 
$\tilde{\mu}$ of $\frk{l}(\frk{\dot{g}},Id,1)''=
\dot{\frk{g}}\otimes\Bb{C}[t,t^{-1}]$ defined by:
\begin{equation}\tilde{\mu}(x\otimes t^j)=(\epsilon_k)^j \mu(x)\otimes t^j,
\qquad {\rm{for}} \ j\in \Bb{Z}, x\in \dot{\frk{g}}. \end{equation}
$\frk{h} = \ocirc{\frk{h}}\otimes 1\oplus\Bb{C}c\oplus\Bb{C}d$ is the standard 
Cartan subalgebra of $\frk{g}$.

\para Let $(\dot{e_i},\dot{f_i})_{i=1,\cdots,n}$ be a system of Chevalley 
generators of $\dot{\frk{g}}$. The simple coroots, ${\dot{\alpha}_i^{\vee}} =
[\dot{e_i},\dot{f_i}]$ for $i=1,\cdots,n$ form a base of $\frk{\dot{h}}$.
Let $\dot{\omega}$ be the Cartan involution of $\dot{\frk{g}}$ given by
$\dot{\omega}(\dot{e_i})=-\dot{f_i}$, $\dot{\omega}\vert_{\dot{\frk{h}}}=-Id$
and  $\dot{\omega}^2=Id$. The simple roots in the base $\dot{\pi}$ of 
$\frk{\dot{g}}$ are enumerated in a manner such that a system of 
representatives of the $\mu$-orbits of $\lb1,\cdots,n\rb$ is $\lb1,\cdots,l
\rb$. Since order of $\mu(\neq Id)$ is $k$(=2 or 3), therefore for any 
$i\in\lb1,\cdots,n\rb$, the cardinality $n_i$ of the $\mu$-orbit of 
${\dot{\alpha}_i^{\vee}}$ is either 1 or $k$. Correspondingly, 
$e_i=\dot{e_i}+\cdots+\mu^{n_i-1}(\dot{e_i})$ except in the case
of $A_{2l}^{(2)}$ where $e_{l}=\sqrt{2}(\dot{e_l}+\dot{e_{l+1}})$; 
$f_i=\ocirc{\omega}(\dot{e_i})$; $\alpha_i^{\vee} = [e_i,f_i]$. Hence,
$(e_i,f_i)_{i=1,\cdots,l}$ is a system of Chevalley generators of the simple 
Lie algebra $\ocirc{\frk{g}}$. Let $\theta_0 \in (\ocirc{\frk{h}})^*$ be the 
highest weight of the irreducible $\ocirc{\frk{g}}$-module $\frk{\dot{g}}_1$. 
Choose $E_0\in (\frk{\dot{g}}_1)_{-\theta_0}$ and put $F_0=-\dot{\omega}(E_0), 
e_0=E_0\otimes t$, $f_0=F_0\otimes t^{-1}$ and $\alpha_0^{\vee}=[e_0,f_0]$. 
Then, $\lb\frk{h}, e_i,f_i, i\in[0,l]\rb$ is a system of generators of 
$\frk{g}(A)$, the Lie algebra associated to a generalized Cartan matrix $A$
(cf. \cite[Theorem 8.3]{Kac}). 
\begin{figure}
{\begin{center}\begin{tabular}{|c|l|c|l|} \hline\hline
Type Aff $k$&\multicolumn{1}{c|}{$\frk{g}$} & 
\multicolumn{1}{c|}{\emph{Dynkin diagram}} &
\multicolumn{1}{c|}{\emph{$\ocirc{\frk{g}}$}}\\ \hline
2&$A_{2}^{(2)}$ & $\underset{\alpha_{0}}{\stackrel{1}{\circ}} \,^{=}_{=}\ >
\underset{\alpha_{1}}{\stackrel{2}{\circ}}$& $A_1$ \\\hline    
2&$ A_{2l}^{\left(2 \right)} \ \left( l \geq 2 \right)$ & $
\underset{\alpha_{0}}{\stackrel{1}{\circ}} \Rightarrow
  \underset{\alpha_{1}}{\stackrel{2}{\circ}}   - \cdots-
  \underset{\alpha_{l-1}}{\stackrel{2}{\circ}} \Rightarrow
  \underset{\alpha_{l}}{\stackrel{2}{\circ}}$ & $B_l$\\\hline
2&$ A_{2l-1}^{\left(2\right)}  \left(l \geq 3 \right)$ &
$  \begin{array}{l} \,\, \qquad \stackrel{1}{\circ} _{\alpha_0}\\
\, \qquad \, \vert \\ 
   \underset{\alpha_{1}}{\stackrel{1}{\circ}} -
     \underset{\alpha{2}}{\circ^{2}} -
   \underset{\alpha_{3}}{\stackrel{2}{\circ}} - \cdots-
   \underset{\alpha_{l-1}}{\stackrel{2}{\circ}} \Leftarrow
   \underset{\alpha_{l}}{\stackrel{1}{\circ}}
 \end{array}$& $C_l$\\\hline
2&$ D_{l+1}^{\left(2\right)} \left(l \geq 2 \right)$&
$ \underset{\alpha_{0}}{\stackrel{1}{\circ}}  \Leftarrow
  \underset{\alpha_{1}} {\stackrel{1}{\circ}}  - \cdots-
  \underset{\alpha_{l-1}}{\stackrel{1}{\circ}}  \Rightarrow
  \underset{\alpha_{l}}{\stackrel{1}{\circ}}$ & $B_l$\\ \hline
2&$E_{6}^{\left(2\right)}$ & 
$\underset{\alpha_{0}}{\stackrel{1}{\circ}}  - 
  \underset{\alpha_{1}}{\stackrel{2}{\circ}}  -
  \underset{\alpha_{2}}{\stackrel{3}{\circ}}  \Leftarrow  
  \underset{\alpha_{3}}{\stackrel{2}{\circ}}  -
  \underset{\alpha_{4}}{\stackrel{1}{\circ}} $& $F_4$\\\hline
3&$D_{4}^{3}$& $\underset{\alpha_{2}}{\stackrel{1}{\circ}}    -
 \underset{\alpha_{1}}{\stackrel{2}{\circ}}   \Lleftarrow
 \underset{\alpha_{2}}{\stackrel{1}{\circ}} $&$G_2$\\
\hline\hline \end{tabular}\end{center}}
\caption{Dynkin diagrams of the twisted affine Kac-Moody Lie algebras.}
\end{figure}

A compact form $\frk{u}(A)$ of $\frk{g}(A)$ is defined as the fixed point set
of $\frk{g}(A)$ under the compact involution $\omega$ defined on
$\frk{\dot{g}}\otimes\Bb{C}[t,t^{-1}]$ as follows:
$\omega( x\otimes t^j) = \dot{\omega}(x)\otimes t^{-j}$.

\para A graph $S(A)$ called the Dynkin diagram of $A$ can be associated to a 
generalized Cartan matrix (GCM) $A$ as explained in \cite{Kac}. A Dynkin 
diagram $S(\ocirc{A})$ is obtained from $S(A)$ by removing the $0^{th}$ vertex.
The corresponding GCM $\ocirc{A}$ is of finite type and $\frk{g}(\ocirc{A})
=\ocirc{\frk{g}}$ is a simple finite dimensional Lie algebra. It is clear that 
$A$ is indecomposable if and only if $S(A)$ is a connected graph. The matrix 
$A$ is determined by the Dynkin diagram and the enumeration of its vertices.
Figure 1 gives the Dynkin diagrams of the twisted affine Kac-Moody Lie 
algebras. The enumeration of the generators is the same as in \cite{Kac} except
in the case $A_{2l}^{(2)}$, where the enumeration is reversed and the case 
$A_2^{(2)}$ where the Dynkin diagram and enumeration has been taken from 
\cite[$\S 3.5$]{Moody}.

\para
Let $\ocirc{\triangle}$ be the root system of $\ocirc{\frk{g}}$. Let 
$\ocirc{\pi}= (\alpha_i)_{i\in[1,l]}$ be a base of 
$(\ocirc{\frk{g}},\ocirc{\frk{h}})$ and $(\ocirc{p_i})_{i\in[1,l]}$ be the 
corresponding dual basis in $\ocirc{\frk{h}}$. Denote by $\ocirc{\triangle}_s$,
$\ocirc{\triangle}_l$ and $\ocirc{\triangle}_+$ the short, long and positive 
roots of $\ocirc{\frk{g}}$. Define an element $\delta\in\frk{h}^*$ by putting 
$\delta(d)=1$ and $\delta(\ocirc{\frk{h}}+\Bb{C}c)=0$. 
Recall from \cite[Proposition 6.3]{Kac} the following description of 
$\triangle^{re}$ 
for $\frk{g}$ not of type $A_{2l}^{(2)}$:
$$ \begin{array}{l}
\triangle^{re} = \lbrace \  \alpha + n \delta \ \mid \  \alpha \in 
\ocirc{\triangle}_s , \ n \in \Bb{Z} \ \rbrace\ \cup   \lbrace \  \alpha + 
nk\delta \ \mid \  \alpha \in \ocirc{\triangle}_l , n \in \mathbb{Z} \rbrace,
\qquad \quad  \text{if}\ k=2\ {\rm{or}}\ 3.\end{array}$$ 
With respect to the enumeration of the simple roots of $A_{2l}^{(2)}$
as given in Figure 1, similar calculations as in \cite[Proposition 6.3]{Kac} 
show that for $\frk{g}$ of type $A_{2l}^{(2)}$, $\triangle^{re} = 
\triangle^{re}_s \cup \triangle^{re}_m\cup \triangle^{re}_{l}$ where,
$$\begin{array}{l}
\triangle^{re}_s = \lbrace \frac{1}{4}(2\alpha + (4n-k)\delta) \mid \ \alpha\in
\ocirc{\triangle}_l, n,k \in \Bb{Z}, 1\leq k\leq 3\rbrace,\\
\triangle^{re}_m = \lbrace \frac{1}{2}(2\alpha + (2n-1)\delta) \mid \ \alpha\in
\ocirc{\triangle}_l, n\in \Bb{Z}\rbrace,\\
\triangle^{re}_l = \lbrace 2\alpha+ n\delta \mid \ \alpha \in 
\ocirc{\triangle}_l, n\in \Bb{Z}\rbrace.\end{array}$$
Let $\triangle^{re}_+ = \lb \alpha \in \triangle^{re}\ \text{with}\ n> 0\rb\cup
\ocirc{\triangle}_+.$ By \cite[$\S$6.4, Proposition 6.4]{Kac}, $a_0\alpha_0 = 
\delta - \theta$ with $\theta\in (\ocirc{\triangle}_+)_s$ for $\frk{g}$ of type
Aff 2 or 3 and not of type $A_{2l}^{(2)}$ and $\alpha_0 = \delta -2\theta$,
with $\theta\in(\ocirc{\triangle_+})_{l}$ for $\frk{g}$  of type
$A_{2l}^{(2)}$. If $a_0,a_1,\ldots,a_l$ are the 
numerical labels of $S(A)$ as in Figure 1, then the element 
$\delta\in \frk{h}^*$ is defined as,
$\delta=\overset{l}{\underset{i=0}{\sum}} a_i\alpha_i$ and we have, 
$$\triangle^{im} =\lb\pm\delta,\pm2\delta,\ldots\rb, \qquad \,
\triangle^{im}_+ =\lb\delta, 2\delta,\ldots \rb.$$ The set $\Pi = 
(\alpha_i)_{i=0,1,\cdots,l}$ is a base of $\triangle.$ \label{p-i}Setting 
$p_0=d$ and $p_i=\ocirc{p_i}+a_id$ for $i=1,\cdots,l$, we obtain a family 
$(p_i)_{i\in[0,l]} \subset \frk{h}$ satisfying $\alpha_j(p_i) =\delta_{i,j}$ 
for $i,j \in [0,l]$. 

Let $\frk{g}_{\gamma}$ be the root space of $\gamma
\in\triangle$. For $\gamma\in\triangle^{re}$, dim $\frk{g}_{\gamma} = 1$.

\label{long/short} For a root $\dot{\alpha}$ of $\frk{\dot{g}}$, if 
$\dot{e}_{\pm \dot{\alpha}}\in \frk{\dot{g}}_{\dot{\alpha}}$ is such that
the $\Bb{C}$-span of $\lb \dot{e}_{\dot{\alpha}},\dot{e}_{-\dot{\alpha}},
H_{\dot{\alpha}}= [\dot{e}_{\dot{\alpha}}, \dot{e}_{-\dot{\alpha}}] \rb$ is 
isomorphic to $\frk{sl}_2$, then given $\alpha+nk\delta\in\triangle^{re}$ with 
$\alpha\in\ocirc{\triangle}_l$, $E_{\pm \alpha, \pm ks} = 
e_{\pm \alpha}\otimes t^{\pm ks}\in\frk{g}_{\pm(\alpha+ks)}$ and 
$H_{\alpha}\in\ocirc{\frk{h}}$ can be suitably chosen such that the 
$\Bb{C}$-span of $\lb E_{\alpha, ks}, H_{\alpha}, E_{-\alpha,-ks}\rb$ is 
isomorphic to $\frk{sl}_2$. If $\alpha\in \ocirc{\triangle}_s$, then there 
exists $\dot{\alpha}\in\dot{\frk{g}}$ such that for $j\in \Bb{Z}_k$ one can 
choose $E_{\pm \alpha, \pm ks+j} = (\eps_k^je_{\pm\mu(\dot{\alpha})}+ \cdots
+\eps_k^{jk}e_{\pm\mu^k(\dot{\alpha})} )\otimes t^{\pm ks+j} \in 
\frk{g}_{\pm(\alpha+ks+j)}$ and $h^j_{\alpha} = \eps_k^jH_{\mu(\dot{\alpha})}+
\cdots+\eps_k^{jk}H_{\mu^k(\dot{\alpha})} \in \frk{\dot{h}}_j$ such that the 
$\Bb{C}$-span of $\lb E_{\alpha,ks+j}, h_{\alpha}^j, E_{-\alpha,-ks-j}\rb$ is 
isomorphic to $\frk{sl}_2$.

\para \label{W(u)} Let $\ocirc{W}$ be the Weyl group of $\ocirc{\frk{g}}$ 
generated by the reflections $(r_i)_{i\in[1,l]}$. Let $T$ be the group of
translations. By \cite[Proposition 6.5]{Kac}, $W = \ocirc{W}\ltimes T$ is the 
Weyl of $\frk{g}$. Since $w(\delta) = \delta$ for all $w \in W$, it follows 
from \cite[$\S$2.7 and Proposition 3.12(b)]{Kac} that given two positive root 
systems $\triangle_+$ and $\triangle_+'$ of a compact form $\frk{u}(A)$ of 
$\frk{g}(A)$, there exists $s\in \ocirc{W}$ such that  
$s.\triangle_+' = \triangle_+$. But for $\alpha \in \ocirc\triangle$, 
the generators $r_{\alpha}$ of $\ocirc{W}$ are interior automorphisms 
(cf. \cite[Lemma 5.3]{M-R}). Hence given two positive root systems 
$\triangle_+$ and $\triangle_+'$ of $\frk{u}(A)$ there exists  
$s\in$ Int$(\frk{u}(A))$ such that $s.\triangle_+' = \triangle_+$.

\section{Automorphisms and Real forms of $\frk{g}$ }

\para Define a group G acting on $\frk{g}$ by the adjoint representation 
Ad:G$\rightarrow$ Aut($\frk{g}$). It is generated by the subgroups 
$U_{\alpha}$, for $\alpha \in \triangle^{re}$, which are isomorphic to the 
additive groups $\frk{g}_{\alpha}$ by an isomorphism $\exp$ such that 
Ad $\circ \ \exp= \exp \circ\ \ad$. 

 A Cartan subalgebra (CSA) of $\frk{g}$ is a maximal ad$_{\frk{g}}$
-diagonalizable Lie subalgebra. The CSA's are all conjugate by G. A 
Borel subalgebra (BSA) is a completely solvable maximal subalgebra of 
$\frk{g}$. $\frk{b}^+=\frk{h}\oplus(\bigoplus_{\alpha\in \triangle^+}
\frk{g}_{\alpha})$ and $\frk{b}^- = \frk{h}\oplus(\bigoplus_{\alpha\in 
\triangle^-}\frk{g}_{\alpha})$ are respectively called the positive and 
negative standard BSA's. The subalgebras $\frk{b}^+$ and $\frk{b}^-$ are not 
conjugate by G. All the BSA's conjugate to $\frk{b}^+$ (respectively 
$\frk{b}^-$) are said to be positive (respectively negative). If $\frk{g}$ is 
indecomposable, all BSAs are either positive or negative.

An automorphism (linear or semi-linear) of $\frk{g}$ acts in a compatible 
manner to Ad on G and hence transforms two conjugate BSAs to two conjugate
BSAs; it is said to be of {\it{first type}} (respectively {\it{second type}})
if it transforms a positive BSA to positive (respectively negative) BSA. If
$\frk{g}$ is indecomposable, all automorphisms are either of first or second 
type. 

\para\label{Aut}Automorphisms of $\frk{g}$ : By \cite{R4}, the group of
automorphisms of $\frk{g}$ is given by 
$$Aut(\frk{g})= [(\lbrace 1, \omega\rbrace \times Aut(A))\ltimes Int(\frk{g})]
\ltimes Tr,$$ where $\omega$ is the Cartan involution of $\frk{g}$, $Aut(A)$
is the group of permutations of [$0,l$] such that 
$a_{\rho i,\rho j} = a_{ij}$ for $i,j \in$ I, Int($\frk{g}$) is the set of 
interior automorphisms of $\frk{g}$ and $Tr$ = Tr($\frk{g},\frk{g}',\frk{c}$) 
is the group of transvections of $\frk{g}$ as defined in \cite[2.4]{R4}.

Let $\frk{h}$ be a standard CSA of $\frk{g}$. A group $\tilde{H}$ is defined 
such that in the complex case, Ad($\tilde{H}$)=$\exp \ad$($\frk{h}$)
(cf.\cite{P-K}) . The group $Int(\frk{g})=Ad(\tilde{H}\ltimes G)$ of interior 
automorphisms of $\frk{g}$ is the image of the semi-direct product of 
$\tilde{H}$ and $G$. Its derived group is the adjoint group $Ad(G)$ (denoted by
$Int(\frk{g'})$). As G acts transitively on the Cartan subalgebra, the group 
$Int(\frk{g'})$ does not depend on the choice of $\frk{h}$. 

\defn Let $Aut_{\Bb{R}}(\frk{g})$ denote the group of automorphisms of 
$\frk{g}$ that are either $\Bb{C}$-linear or semilinear(i.e., $\phi(\lambda x)=
\bar{\lambda} \phi(x), \ \forall\ \lambda \in \Bb{C}, x \in \frk{g}$). 
$Aut(\frk{g})$ is an index 2 normal subgroup of $Aut_{\Bb{R}}(\frk{g})$.

A semi-involution of $\frk{g}$ is a semi-linear automorphism of order 2. For 
all semi-involutions $\sigma'$ we have a decomposition, 
$Aut_{\Bb{R}}(\frk{g}) = \lbrace 1, \sigma'\rbrace \ltimes Aut(\frk{g})$.
If $\sigma'$ is a semi-involution of $\frk{g}$, the real Lie algebra 
$\frk{g}_{\Bb{R}} = \frk{g}^{\sigma'}$ is a real form of $\frk{g}$, in the 
sense that there exists an isomorphism of the complex Lie algebras 
$\frk{g}_{\Bb{R}}\otimes_{\Bb{C}} \Bb{C}$ and $\frk{g}$; further, $\sigma'$ is
the conjugation of $\frk{g}$ with respect to $\frk{g}_{\Bb{R}}$. Thus there 
exists a bijective correspondence between the semi-involutions and real forms.
The {\it{standard normal(or split)}} real form of $\frk{g}$ is the real 
Lie algebra generated by $e_i,f_i,\alpha_i^{\vee}$ and $d$. The corresponding
semi-involution $\sigma_n'$ is called the {\it{normal semi-involution}}.
Note that $\sigma_n'$ is the restriction of $\dot{\sigma}_{n}'\otimes conj$
on $\frk{g}''=\big(\dot{\frk{g}}\otimes\Bb{C}[t,t^{-1}]\big)^{\tilde{\mu}}$, 
where $\dot{\sigma}_n'$ is the normal semi-involution of $\frk{\dot{g}}$ and
$conj(P(t))=\bar{P}(t)$. $\sigma_n'$ commutes with the standard Cartan 
involution $\omega$.

The {\it{standard Cartan semi-involution}} $\omega'$ of $\frk{g}$ is the unique
semi-involution of $\frk{g}$ such that $\omega'(e_i) = - f_i$, and $\omega'(d)$
$= -d$. Hence $\omega' = \sigma_n'\omega = \omega\sigma_n'$. In the standard 
realization of $\frk{g}$, $\omega'$ induces on $\frk{g}''$ the restriction of 
$\dot{\omega}\otimes\iota'$, where $\dot{\omega}$ is the Cartan 
semi-involution of $\dot{\frk{g}}$ and $\iota'(P(t))=\bar{P}(t^{-1})$. 
All conjugates of $\omega'$ are called {\it{Cartan semi-involutions (CSI)}} or 
{\it{compact semi-involutions}}; these are semi-involutions of the second type.
The corresponding real forms are called the {\it{compact real forms}}. It is 
clear that, for all affine Kac-Moody Lie algebras, there exists, upto a 
conjugation, a unique compact real form.

\para{\it{Real Forms of $\frk{g}$}} : 
The real form corresponding to a semi-involution of first type (SI1) 
(respectively of second type (SI2)) is said to be \emph{almost split} 
(respectively \emph{almost compact}) real form. Upto a conjugation, a 
classification of the almost split real forms was given in \cite{1type} and a 
classification of the almost compact real forms was given in \cite{M-R}. 

\para{\it{Cartan subalgebra of a real form of $\frk{g}$}}: Let 
$\frk{g}_{\Bb{R}}$ be a real form of the complex Lie algebra $\frk{g}$. A Lie 
subalgebra $\frk{h}_0$ of $\frk{g}_{\Bb{R}}$ is called the 
{\emph{Cartan subalgebra of $\frk{g}_{\Bb{R}}$}} if the complexification, 
$\frk{h}_0\otimes \Bb{C}$ is a Cartan subalgebra of $\frk{g}$.

\para \emph{Cartan Involutions}:
Let $\sigma'$ be a SI2 of $\frk{g}$, and let $\frk{g}_{\Bb{R}} = 
\frk{g}^{\sigma'}$ be the corresponding almost compact real form.
 A CSI $\vartheta$ that commutes with $\sigma'$ is said to be \emph{adapted}
to $\sigma'$ or $\frk{g}_{\Bb{R}}$. The involution $\sigma = \sigma'\vartheta$
(respectively its restriction $\vartheta_{\Bb{R}}$ to  $\frk{g}_{\Bb{R}}$ or  
$\sigma'_{\Bb{R}}$ to  $\frk{u} = \frk{g}^{\vartheta}$ ) is said to be the 
Cartan involution of $\sigma'$ (respectively of $\frk{g}_{\Bb{R}}$ or of 
$\frk{u}$). The algebra of fixed points $\frk{k} = \frk{g}_{\Bb{R}}^{\sigma} = 
\frk{g}_{\Bb{R}}\cap\frk{u} = \frk{u}^{\sigma}$ is the \emph{maximal compact 
subalgebra} of $\frk{g}_{\Bb{R}}$. $\frk{g}_{\Bb{R}} = \frk{k} \oplus \frk{p}$;
$\frk{u} = \frk{k} \oplus i\frk{p}$ is a Cartan decomposition of
$\frk{g}_{\Bb{R}}$ and $\frk{u}$ into eigenspaces of $\sigma$. A Cartan 
subalgebra $\frk{h}$ of $\frk{g}$ is said to be maximally compact for 
$\sigma'$ (or $\frk{g}_{\Bb{R}}$) if it is stable under $\sigma'$ and if 
$-\sigma'$ stabilizes a base of $\triangle(\frk{g},\frk{h})$. 

\begin{prop}\label{CSI}\cite[Proposition 2.8]{R4} Let $\sigma'$ be a SI2 of 
$\frk{g}$ and let $\frk{g}_{\Bb{R}} = \frk{g}^{\sigma'}$.\\
1. There exists CSI's $\vartheta$ adapted to $\sigma'$ and
maximally compact CSA's $\frk{h}$ for $\sigma'$.\\
2. For all maximally compact CSA's $\frk{h}$ for $\sigma'$, 
there exists a CSI $\vartheta$ adapted to $\sigma'$ that stabilizes $\frk{h}$,
and it is unique upto an interior automorphism fixing $\frk{h}$ and commuting 
with $\sigma'$.\\ 
3. For all CSI $\vartheta$ adapted to $\sigma'$, there exists a maximally 
compact Cartan subalgebra for $\sigma'$ stabilized by $\vartheta$, and it is 
unique upto an interior automorphism commuting with $\sigma'$ and $\vartheta$.
\end{prop}

\begin{prop}\label{1-2pair}\cite[Proposition 2.9]{R4} Consider:\\
1. $\sigma$, an involution of first type of $\frk{g}$.\\
2. the pairs ($\sigma',\frk{h}$) formed of a semi-involution of second type 
$\sigma'$, which is not Cartan, and a maximally compact Cartan sub-algebra for 
$\sigma'$.\\ 
3. the relation  ($\sigma',\frk{h}$) $\sim \sigma$ \iff $\sigma$ commutes 
with $\sigma'$, stabilizes $\frk{h}$ and is such that $\sigma\sigma'$ is a 
Cartan semi-involution of $\frk{g}$.\\
This relation induces a bijection between the conjugacy classes (under 
$Int_{Tr}(\frk{g})$ or $Aut(\frk{g})$) of involutions of the first type, 
$\sigma$ and the pairs  ($\sigma',\frk{h}$).\\
{\bf{Note}}: With notations as above, if  ($\sigma',\frk{h}$) $\sim \sigma$, 
then we say that $\sigma$ is adapted to $\sigma'$ or $\frk{g}_{\Bb{R}}=
\frk{g}^{\sigma'}$.
\end{prop}
\notn In what follows, $\dot{\frk{s}}$ will denote a finite dimensional
semisimple Lie algebra over $\Bb{C}$.

\para {\it{Adapted Realization of $\frk{g}$}}: Recall the following 
definitions from \cite{M-R}:
\begin{defn} Let $\sigma$ be finite order automorphism of $\frk{g}$ of first 
type and $\frk{h}$ a maximally fixed Cartan subalgebra for $\sigma$. A 
realization of $\frk{g}$ on which $\sigma$ preserves the $\Bb{Z}$-gradation 
and for which $\frk{h}$ is the standard Cartan subalgebra is said to be 
{\it{ almost adapted to ($\sigma, \frk{h}$)}}. An almost adapted realization 
$\frk{l}(\frk{\dot{s}},\zeta,\eps_m)$ for ($\sigma, \frk{h}$) on which $\sigma$
commutes with the translation map $\mathcal{T}:x\mapsto x\otimes t^m$, is said 
to be {\it{adapted to $\sigma$}} (respectively to ($\sigma, \frk{h}$)).
\end{defn}

\begin{defn} Let $\sigma'$ be finite order automorphism of $\frk{g}$ of second 
type, $\frk{h}$ a maximally compact Cartan subalgebra for $\sigma'$ and $\sigma
$ an involution of first kind associated to the pair ($\sigma',\frk{h}$). A 
realization $\frk{l}(\frk{\dot{s},\mu,\varepsilon_m})$ of $\frk{g}$ which is
adapted to ($\sigma,\frk{h}$) is said to be {\it{adapted to $(\sigma',\sigma,
\frk{h})$}} if there exists an involution $\dot{\sigma}$ and a semi-involution 
$\dot{\sigma}'$ of $\dot{\frk{s}}$ commuting with $\mu$ such that 
$\dot{\sigma}\dot{\sigma}'$ is a Cartan semi-involution of $\frk{\dot{s}}$, 
$\sigma = \dot{\sigma} \otimes 1$ and $\sigma' = \dot{\sigma}'\otimes \iota'$ 
on the realization $l''(\frk{\dot{s}}, \mu,\varepsilon_m)$; and finally 
$\sigma'(c)=-c$ and $\sigma'(d) = -d.$ 
\end{defn}

It is known from \cite[Proposition 3.4, Theorem 3.5 and Proposition 3.9]{M-R}
that given a finite order automorphism $\sigma$ of first type (respectively 
$\sigma'$ of second type) of $\frk{g}$, there exists realizations of $\frk{g}$ 
adapted to $\sigma$ (respectively adapted to $\sigma'$).

\para\label{sigma} 
Let $\sigma$ be an involution adapted to a semi-involution of second type
$\sigma'$. Since $\sigma$ is an involution of $\frk{g}$ of the first kind, by 
\cite[Chapter II]{Bausch} upto an interior automorphism, either 
$\sigma$ = $\rho H$ or $\sigma = H$, where $\rho$ is a diagram automorphism of 
$\frk{g}=\frk{l}(\dot{\frk{g}},\mu,\varepsilon_k)$ and $H$ is an interior 
automorphism of $\frk{g}$ of the form $\exp i \pi \ad(h_0),$ with 
$h_0 \in \frk{h}^{\rho}_{\Bb{Z}}$ where $\frk{h}^{\rho}_{\Bb{Z}} = 
\lbrace x \in \frk{h}^{\rho} \mid \alpha(x) \in \Bb{Z}, \ \forall 
\ \alpha \in \triangle \rbrace$. Further $h \in \frk{h}_{\Bb{Z}}$ can be 
written as $h = h' + \eta d,$ with $h'\in (\dot{\frk{g}})^{\mu}$ and $\eta \in 
\Bb{Z}$.

\begin{prop}\label{H}\cite[Proposition 2.12]{M-R} With notations as above:\\
1. If $\eta$ is even or $k=2$ (i.e., $\frk{g}$ is of type Aff 2) then the 
interior involution $H$ respects the $\Bb{Z}-$gradation of $\frk{g} = 
l(\frk{\dot{g}},\mu, \varepsilon_k)$, commutes with the map 
$\mathcal{T}:x\mapsto x\otimes t^{k}$ on $\frk{g}''$ and there exists an 
interior involution $\dot{H}$ of $\frk{\dot{g}}$ commuting with $\mu$ such that
$H = (\mu^{\eta}\dot{H})\otimes 1$ on $l''(\frk{\dot{g}},\mu,\varepsilon_k) = 
\frk{g}''$.\\ 
2. If $\eta$ is odd and $k\neq 2$, then the interior involution $H$ induces
on $\frk{g}''$ the automorphism $t^k\mapsto -t^k$ and $H$ in this case acts
on the adapted realization 
$\frk{l}(\frk{\dot{g}}\times\frk{\dot{g}},\zeta,\varepsilon_{2k})$ of 
$\frk{g}$, with $\zeta(x,y)=(y,\mu{(x)})$.\\ 
3. Each diagram automorphism of a twisted affine Kac-Moody Lie algebra 
$\frk{g}=\frk{l}(\frk{\dot{g}},\mu,\varepsilon_m)$ is $\Bb{A}-$linear, where 
$\Bb{A}$ is the algebra generated by the translation maps 
$\mathcal{T}^{\pm}:x\mapsto x\otimes t^{\pm m}$. 
\end{prop}

\noindent The following are some examples of realizations of $\frk{g} = 
\frk{l}(\frk{\dot{g}},\mu,\varepsilon_k)$, $k=2,3$, adapted to 
($\sigma,\sigma',\frk{h}$).\\
Example 1. Let  $\frk{g} = \frk{l}(\frk{\dot{g}},\mu,-1)$ where
$\dot{\frk{g}}$ is of type $A_{2l-1}$.\\
i. Let $\sigma = H$, where $H$ is an interior automorphism of $\frk{g}$. 
Then by Proposition \ref{H}(1), $\frk{l}(\frk{\dot{g}},\mu,-1)$ is adapted to 
($\sigma,\sigma',\frk{h}$).\\
ii. Let $\sigma = \rho$, where $\rho$ is a diagram automorphism of 
$A_{2l-1}^{(2)}$ such that $\rho(\alpha_0)=\alpha_1$, $\rho(\alpha_1)=\alpha_0$
and $\rho(\alpha_i)=\alpha_i$ for $i\in[2,l]$. In the given realization 
since $\sigma(E_0\otimes t)=\rho(E_0\otimes t) = e_1\otimes 1$, $\sigma$ 
clearly cannot be written in the form $\dot{\sigma}\otimes 1$ for any
 $\dot{\sigma}\in Aut(\frk{\dot{g}})$. Hence using \cite[1.7]{M-R}, with 
$H_1 =\exp \frac{i\pi}{2}\ad(\ocirc{p_1})$ we consider the realization 
$\frk{l}(\frk{\dot{g}},\mu H_1,\varepsilon_4)$ of $\frk{g}$ which is 
isomorphic to $\frk{l}(\frk{\dot{g}},\mu,-1)$ via the map 
$\psi: \frk{l}(\frk{\dot{g}},\mu,-1) \mapsto 
\frk{l}(\frk{\dot{g}},\mu H_1,\varepsilon_4) $ defined as follows:
$$\begin{array}{c} 
\psi(x\otimes t^j)=x\otimes t^{2j+N}+\delta_{j,0} (\ocirc{p_1},x) C; \quad
\text{for }\ x\in\frk{\dot{g}}_j\ \text{such that}\ [\ocirc{p_1},x]=Nx,\\
\psi(c)=2C, \qquad \psi(d)=(D-\ocirc{p_1})/2,
\end{array}$$where $C$ and $D$ respectively denote the central and gradation 
elements of $\frk{l}(\frk{\dot{g}},\mu H_1,\varepsilon_4)$. Since
[$\ocirc{p_1},E_0$]=$-E_0$ and [$\ocirc{p_1},e_1$]=$e_1$, it is easy to see 
that under the new realization $e_0=E_0\otimes t$ is mapped to $E_0\otimes t$
and $e_1=e_1\otimes1$ is mapped to $e_1\otimes t$ and hence 
$\rho$ can be written in the form $\dot{\sigma}\otimes1$ where $\dot{\sigma}\in
Aut(\frk{\dot{g}})$ is an involution such that $\dot{\sigma}(E_0)=e_1$, and 
$\dot\sigma(e_i)=e_i$ for $i\in[2,l]$. Taking 
$\dot{\sigma}'=\dot{\sigma}\dot{\omega}$ it can be easily
seen that $\dot{\sigma} \mu H_1 =\mu H_1 \dot\sigma$ on $\frk{\dot{g}}$. Hence
the realization $\frk{l}(\frk{\dot{g}},\mu H_1,\varepsilon_4)$ is adapted
to ($\rho,\rho\omega',\frk{h}$). 

\noindent Example 2.  
Let  $\frk{g} = \frk{l}(\frk{\dot{g}},\mu,\varepsilon_3)$ where $\dot{\frk{g}}$
is of type $D_4$.\\
i. Let $\sigma = \exp i\pi\ad(p_1)$ or $\exp i\pi\ad(\ocirc{p_2})$.
Then by Proposition \ref{H}(1), $\frk{l}(\frk{\dot{g}},\mu,\varepsilon_3)$ is 
adapted to ($\sigma,\sigma',\frk{h}$).\\ 
ii. Let $\sigma = \exp i\pi\ad(p_0)$. Then by Proposition \ref{H}(2), 
$\frk{l}(\frk{\dot{g}}\times\frk{\dot{g}},\zeta,\varepsilon_{6})$ is an 
adapted realization of $\frk{l}(\frk{\dot{g}},\mu,\varepsilon_{3})$. 
It is known from \cite{M-R} that for $\zeta\in 
Aut(\frk{\dot{g}}\times\frk{\dot{g}} )$ as defined above, 
$\phi_2:\frk{l}(\frk{\dot{g}},\mu,\varepsilon_{3})\mapsto 
\frk{l}(\frk{\dot{g}}\times\frk{\dot{g}},\zeta,\varepsilon_{6})$ is defined by:
$$\begin{array}{c} \phi_2(x\otimes t^j)=(x,\varepsilon_{6}x)\otimes t^j,\qquad
\phi_2(c)=2c_2, \qquad \phi_2(d)=d_2, \end{array}$$ gives an isomorphism 
between the two realizations of $\frk{g}$. Here $c_2$ and $d_2$ respectively 
denote the central and gradation elements of 
$\frk{l}(\frk{\dot{g}}\times\frk{\dot{g}},\zeta,\varepsilon_6)$.

\para\label{Bw} Let $\sigma'$ be a semi-involution of second type of $\frk{g}$ 
and let $\frk{g}_{\Bb{R}} = (\frk{g})^{\sigma'}$ be the corresponding almost 
compact real form of $\frk{g}$. Let $\omega'$ is a Cartan semi-involution of 
$\frk{g}$ adapted to $\sigma'$. Then $B_{\omega'}$ defined by $B_{\omega'}(X,Y)
= -B(X,\omega'(Y))$ is a positive definite hermitian form on $\frk{g}_{\Bb{R}}$
and given $X,Y,Z \in \frk{g}_{\Bb{R}}$, we have:
\begin{equation}\begin{array}{cl}
B_{\omega'}(\ad\ \sigma X(Y),Z) &= - B([\sigma X,Y],\omega' Z)\\
&=\ \ B(Y, [\sigma X,\omega' Z])\\ &=\ \ B(Y, \omega'[\sigma' X,Z])\\
&= - B_{\omega'}(Y,\ad\ \sigma' X(Z)).
\end{array}\end{equation} Since for all $X \in \frk{g}_{\Bb{R}} = 
\frk{g}^{\sigma'}$, $\sigma' X =X$, therefore with respect to $B_{\omega'}$, we
get $$(\ad\ \sigma X)^{*} = -\ad\ \sigma' X  = -\ad X\ \ \forall \ X \in 
\frk{g}_{\Bb{R}} = \frk{g}^{\sigma'}.$$
This implies that for $X \in \frk{k} = \frk{g}^{\sigma}$, $\ad X$ is 
skew-symmetric with purely imaginary eigenvalues and  for $X \in \frk{p}$, 
$\ad X$ is symmetric with real eigenvalues.

Let $\frk{l}(\frk{\dot{s}},\mu,\varepsilon_m)$ be a 
realization of $\frk{g}$ adapted to ($\sigma',\sigma,\frk{h}$) where 
$\frk{\dot{s}}$ is a complex semi-simple Lie algebra. Let 
$\dot{\sigma},\ \dot{\sigma}'\in Aut(\dot{\frk{s}})$ such that $\sigma = 
\dot{\sigma}\otimes 1$ and $\sigma' = \dot{\sigma}'\otimes \iota'$. Let 
$\frk{\dot{s}}_{\Bb{R}}$ = $\dot{\frk{s}}^{\dot{\sigma}'}$ and 
$\frk{\dot{s}}_{\Bb{R}} = \frk{\dot{k}}\oplus \frk{\dot{p}}$ be its Cartan 
decomposition with respect to $\dot{\sigma}$.

\begin{lem}\label{hR}With the above notations, let $\frk{l}(\frk{\dot{s}},\mu,
\varepsilon_m)$ be a realization of $\frk{g}$ adapted to 
($\sigma',\sigma,\frk{h}$) and $\frk{\dot{t}}$ a $\mu$-stable maximal 
abelian subspace of $\frk{\dot{k}} = (\frk{\dot{s}}_{\Bb{R}})^{\dot{\sigma}}$. 
Then $Z_{\frk{\dot{s}}_{\Bb{R}}}(\frk{\dot{t}})$ is a $\dot\sigma$ stable 
subalgebra of $\frk{\dot{s}}_{\Bb{R}}$ of the form  $Z_{\frk{\dot{s}}_{\Bb{R}}}
(\frk{\dot{t}}) = \frk{\dot{t}} \oplus \frk{\dot{a}}$, where $\frk{\dot{a}}
\subset \frk{\dot{p}}$ and $\frk{h}_{\Bb{R}}'' = (\frk{\dot{t}})^{\mu}\otimes1
\oplus(\frk{\dot{a}})^{\mu}\otimes1$ is a Cartan subalgebra of 
$\frk{g}_{\Bb{R}}''$. \end{lem} 
\pf Given a finite-dimensional semisimple Lie algebra 
 $\frk{\dot{s}}$ over $\Bb{C}$, by \cite[Proposition 6.60]{Knapp}, 
$Z_{\frk{\dot{s}}_{\Bb{R}}}(\frk{\dot{t}})$ is a $\dot{\sigma}$-stable Cartan 
subalgebra of $\frk{\dot{s}}_{\Bb{R}}$ of the form $Z_{\frk{\dot{s}}_{\Bb{R}}}
(\frk{\dot{t}}) = \frk{\dot{t}} \oplus \frk{\dot{a}}$, with $\frk{\dot{a}}
\subset \frk{\dot{p}}$. The complex Lie algebra 
$\mathcal{Z} = \big(Z_{\frk{\dot{s}}_{\Bb{R}}}(
\frk{\dot{t}})\big)_{\Bb{C}}$ being a Cartan subalgebra of $\frk{\dot{s}}$, is
a $\mu-$stable reductive subalgebra of $\frk{\dot{s}}$. 
Consider the infinite abelian subalgebra
$\frk{z}''$=$\frk{l}''(\mathcal{Z},\mu,\varepsilon_m)$ of $\frk{g}''$. 
If $x\in\mathcal{Z}\subset\dot{\frk{s}}$,
$$\begin{array}{l}
\sigma(x\otimes t^s) = (\dot{\sigma}\otimes1)(x\otimes t^s)
= \dot{\sigma}(x)\otimes t^s. \end{array}$$ As $Z_{\frk{\dot{s}}_{\Bb{R}}}
(\frk{\dot{t}})$ is stable under the linear involution $\dot{\sigma}$,
$\dot{\sigma}(x)\in \mathcal{Z}$ for $x\in\mathcal{Z}$. Hence $\frk{z}''$ is
$\sigma$-stable and consequently we have $\frk{z}'' = \big(\frk{z}''\cap
\frk{l}''(\frk{\dot{p}}_{\Bb{C}},\mu,\varepsilon_m)\big)\oplus
\big(\frk{z}''\cap\frk{l}''(\frk{\dot{k}}_{\Bb{C}},\mu,\varepsilon_m)\big)$.
But by \cite[Proposition 5.1]{M-R}, the semisimple elements of 
$\frk{z}''\cap\frk{l}''(\frk{\dot{p}}_{\Bb{C}},\mu,\varepsilon_m)$ are 
contained in $(\frk{\dot{a}})^{\mu}\otimes1$ and the semisimple elements of 
$\frk{z}''\cap\frk{l}''(\frk{\dot{k}}_{\Bb{C}},\mu,\varepsilon_m)$ are 
contained in $(\frk{\dot{t}})^{\mu}\otimes1$. Hence an element $x \in 
\frk{z}''\cap\frk{g}_{\Bb{R}}''$ is semisimple if and only if $x \in (\frk{\dot
{t}})^{\mu}\otimes1\oplus(\frk{\dot{a}})^{\mu}\otimes1$. Thus
$(\frk{\dot{t}})^{\mu}\otimes1\oplus(\frk{\dot{a}})^{\mu}\otimes1$ is a 
semi-simple abelian subalgebra of $\frk{g}_{\Bb{R}}''$. The lemma will now 
follow if we prove that $ \big((\frk{\dot
{t}})^{\mu}\otimes1\oplus(\frk{\dot{a}})^{\mu}\otimes1\big)_{\Bb{C}} =
\big((\frk{\dot{t}}\oplus\frk{\dot{a}})^{\mu}\otimes1\big)_{\Bb{C}}$ is a 
Cartan subalgebra of $\frk{g}''$. 

From the definition of the action of $\mu$ on the Lie algebra 
$\frk{l}(\frk{\dot{s}}, Id, 1)''$ it clearly follows that,
$$\big((\frk{\dot{t}}\oplus\frk{\dot{a}})^{\mu}\otimes1\big)_{\Bb{C}}
= (\frk{\dot{t}}\oplus\frk{\dot{a}})^{\mu}\otimes \Bb{C} =
\big((\frk{\dot{t}}\oplus\frk{\dot{a}}) \otimes {\Bb{C}}\big)^{\mu}=
\big(Z_{\frk{\dot{s}}_{\Bb{R}}}(\frk{\dot{t}})\otimes \Bb{C}\big)^{\mu} =
\mathcal{Z}^{\mu}.$$ But $\mathcal{Z}$ is a Cartan subalgebra of 
$\frk{\dot{s}}$ and it is known from \cite[Chapter 8]{Kac}, that
$\mathcal{Z}^{\mu}$ is a Cartan subalgebra of $\frk{g}''= 
\frk{l}(\frk{\dot{s}},\mu,\epsilon_{m})''$, whenever $\mathcal{Z}$ is a Cartan 
subalgebra of $\frk{\dot{s}}$. Hence the claim.
$\hfill \blacksquare$ \endpf

\para Let $\frk{g}_{\Bb{R}}$ be an almost compact noncompact real form of 
$\frk{g}$ corresponding to the SI2 $\sigma'$, 
$\frk{l}(\frk{\dot{s}},\mu,\varepsilon_m)$ a realization of $\frk{g}$ adapted 
to $\sigma'$, $\sigma$ a Cartan involution associated to $\sigma'$ and 
$\frk{g}_{\Bb{R}} = \frk{k}\oplus\frk{p}$ the Cartan decomposition of 
$\frk{g}_{\Bb{R}}$ with respect to $\sigma$. Let $\frk{h}_{\Bb{R}}'' = 
(\frk{\dot{t}})^{\mu}\otimes1\oplus(\frk{\dot{a}})^{
\mu}\otimes1$, with $(\frk{\dot{t}})^{\mu} \subset \frk{\dot{k}}$ and  
$(\frk{\dot{a}})^{\mu} \subset \frk{\dot{p}}$ be a $\sigma-$stable Cartan 
subalgebra  of an almost compact real form $\frk{g}_{\Bb{R}}'' = 
\frk{l}(\frk{\dot{s}}_{\Bb{R}},\mu,\varepsilon_m)''$.  $\frk{h}_{\Bb{R}}''$
is said to be maximally compact if dim $(\frk{\dot{t}})^{\mu}$ is as large as 
possible. It clearly follows from paragraph {\bf{\ref{Bw}}} and Lemma \ref{hR},
that all the real roots of a maximally compact Cartan subalgebra 
$\frk{h}_{\Bb{R}}$ of an almost compact real form $\frk{g}_{\Bb{R}}$ have real 
eigenvalues on $(\frk{\dot{a}})^{\mu}\otimes1$ and imaginary eigenvalues on 
$(\frk{\dot{t}})^{\mu}\otimes1$. A real root is called {\it{$\sigma$-real}} if 
it takes real values on $\frk{h}_{\Bb{R}}$(i.e., vanishes on 
$(\frk{\dot{t}})^{\mu}\otimes1$), {\it{$\sigma$-imaginary}} if it takes purely 
imaginary values on $\frk{h}_{\Bb{R}}$ (i.e., vanishes on 
$(\frk{\dot{a}})^{\mu}\otimes1$) and {\it{$\sigma$-complex}} \ otherwise. For 
any $\alpha \in \triangle^{re}$, let $X_{\alpha} \in \frk{g}_{\alpha}$. Then 
$$ [H, \sigma X_{\alpha}] = \sigma[\sigma^{-1}H,X_{\alpha}] = 
\alpha(\sigma^{-1}H) \sigma X_{\alpha}.$$ Hence $\sigma\alpha(H) =\alpha(
\sigma^{-1}H) $ is a root. If $\alpha$ is $\sigma$-imaginary, then 
$\sigma\alpha = \alpha$. In this case, $\frk{g}_{\alpha}$ is $\sigma$-stable, 
and we have $\frk{g}_{\alpha} = (\frk{g}_{\alpha}\cap \frk{k})\oplus
(\frk{g}_{\alpha}\cap \frk{p})$. Since $\frk{g}_{\alpha}$ is 1-dimensional, 
$\frk{g}_{\alpha}\subset\frk{k}$ or $\frk{g}_{\alpha}\subset \frk{p}$. A
$\sigma$-imaginary real root $\alpha$ is said to be {\it{compact}} if 
$\frk{g}_{\alpha}\subset \frk{k}$, {\it{noncompact}} if 
$\frk{g}_{\alpha}\subset\frk{p}$. 

\para Let $\alpha\in\triangle^{re}$ be $\sigma$-real. Since $\sigma\alpha = 
-\alpha$, $\omega'(\alpha)=-\alpha$ and $\sigma\sigma'$ is a 
Int($\frk{g}$)-conjugate of $\omega'$, therefore 
$\sigma'H_{\alpha} = H_{\alpha}$.
Consequently $H_{\alpha}\in \frk{p}\subset\frk{g}_{\Bb{R}}$ is $\ad_{\frk{g}}-$
diagonalizable with real eigenvalues. For $X_{\pm \alpha} \in 
\frk{g}_{\pm \alpha}$ and $H\in \frk{h}^{\sigma'}$, we have, 
$$[H,\sigma'X_{\pm \alpha}] = \sigma'[\sigma'(H),X_{\pm \alpha}] = 
\pm \alpha(\sigma'(H))\sigma'X_{\pm\alpha} = \pm \alpha(H)\sigma'X_{\pm\alpha}.
$$ Therefore $\sigma'(\frk{g}_{\pm\alpha})\subset
\frk{g}_{\pm\alpha}$. Fixing $E_{\alpha}\in \frk{g}_{\alpha}$, $E_{-\alpha} \in
\frk{g}_{-\alpha}$ so that $B(E_{\alpha},E_{-\alpha}) =1$, we get $[E_{\alpha},
E_{-\alpha}]= H_{\alpha}$. Since for $\alpha\in\triangle^{re}$, 
dim $\frk{g}_{\pm\alpha}$=1, the fact that $\sigma'$ is an involution
implies that $\sigma'(E_{\pm \alpha}) = E_{\pm \alpha}$ or 
$\sigma'(E_{\pm \alpha}) = -E_{\pm \alpha}$.

\noindent {\it{Case 1}}: If $\sigma'(E_{\pm \alpha}) = E_{\pm \alpha}$, then 
$\Bb{R}E_{\alpha}\oplus \Bb{R}H_{\alpha}\oplus\Bb{R}E_{-\alpha} \subset 
\frk{g}_{\Bb{R}}.$ By suitably scaling $E_{\alpha}$ and $E_{-\alpha}$ we can 
find $E_{\alpha} \in \frk{g}_{\alpha}$, $E_{-\alpha} \in \frk{g}_{-\alpha}$ and
$H_{\alpha} \in \frk{h}^{\sigma'}$ such that 
$$\begin{array}{ccc}[H_{\alpha},E_{\alpha}] = 
2E_{\alpha}, & [H_{\alpha},E_{-\alpha}]= - 2E_{-\alpha}, & 
[E_{\alpha},E_{-\alpha}]=H_{\alpha}.\end{array}$$ 
From the definition of the Cartan semi-involution and the Cartan involution 
$\sigma$ adapted to $\sigma'$ it now clear that $\sigma(E_{\alpha})=
-E_{-\alpha}$ and $\sigma(E_{-\alpha})=-E_{\alpha}$. Hence 
$E_{\alpha}-E_{-\alpha} \in \frk{k}$. Corresponding to such a $\sigma$-real 
root $\alpha$ define an automorphism $D^{r}_{\alpha}$ as follows: 
$$D^r_{\alpha} = Ad(\exp\frac{\pi}{4}(-E_{-\alpha}-E_{\alpha})).$$ Following 
similar calculations as done in \cite[Proposition 6.52]{Knapp}, we get 
$D^r_{\alpha}(H_{\alpha}) = (E_{\alpha} - E_{-\alpha})$. 

\noindent {\it{Case 2}}: 
If $\sigma'(E_{\pm \alpha}) = -E_{\pm \alpha}$, then $\Bb{R}iE_{\alpha}\oplus 
\Bb{R}H_{\alpha}\oplus\Bb{R}iE_{-\alpha} \subset \frk{g}_{\Bb{R}}.$ 
By suitably scaling $E_{\alpha}$ and $E_{-\alpha}$ we can find $E_{\alpha} \in 
\frk{g}_{\alpha}$, $E_{-\alpha} \in \frk{g}_{-\alpha}$ and $H_{\alpha} \in 
\frk{h}^{\sigma'}$ such that $$\begin{array}{ccc}[H_{\alpha},iE_{\alpha}] = 
2iE_{\alpha}, & [H_{\alpha},iE_{-\alpha}]= - 2iE_{-\alpha}, & 
[iE_{\alpha},iE_{-\alpha}]=-H_{\alpha}.\end{array}$$ 
From the definition of the Cartan semi-involution and the Cartan involution 
$\sigma$ adapted to $\sigma'$ it now clear that $\sigma(E_{\alpha})=E_{-\alpha}
$ and $\sigma(E_{-\alpha})=E_{\alpha}$. Hence $i(E_{\alpha}+E_{-\alpha}) \in 
\frk{k}$. Corresponding to such a $\sigma$-real root $\alpha$ define an 
automorphism $D^{im}_{\alpha}$ as follows: $$D^{im}_{\alpha} = 
Ad(\exp i\frac{\pi}{4}(E_{-\alpha}-E_{\alpha})).$$ Similar calculations as in 
\cite[Proposition 6.52]{Knapp}, show that $D^{im}_{\alpha}(H_{\alpha}) = 
i(E_{\alpha}+E_{-\alpha})$. 

Given a $\sigma$-real root $\alpha$, $D^r_{\alpha}(H_{\alpha})$ and 
$D^{im}_{\alpha}(H_{\alpha})$ being the image of semisimple elements are 
semisimple. Hence, from above discussion it follows that, 
$$ \begin{array}{l}
\bullet\ \text{if  $\sigma'$ fixes $\frk{g}_{\alpha}\oplus\frk{g}_{-\alpha}$, 
then $D^r_{\alpha}$ increases the dimension of $(\frk{\dot{t}})^{\mu}\otimes1$ 
by 1;}\\
\bullet\ \text{if  $-\sigma'$ fixes $\frk{g}_{\alpha}\oplus\frk{g}_{-\alpha}$, 
then $D^{im}_{\alpha}$ increases the dimension of 
$(\frk{\dot{t}})^{\mu}\otimes1$ by 1}.\end{array}$$
Now replacing $d_{\alpha}$ by $D^r_{\alpha}$ and $D^{im}_{\alpha}$ 
appropriately, the same proof as \cite[Proposition 3.8]{B1} shows :
\begin{lem}\label{M.c} A $\sigma$-stable Cartan subalgebra $\frk{h}_{\Bb{R}}''$
of $\frk{g}_{\Bb{R}}''$ is maximally compact if and only if there are no 
$\sigma$-real real roots in $\frk{h}_{\Bb{R}}''$.\footnote{An independent proof of Lemma 3.15 is given in \cite[2.6 iii]{M-R2}}\end{lem}
\begin{rem}\label{h.decomp} As a consequence, we have $\frk{h}_{\Bb{R}}$=
$\frk{h}\cap\frk{g}_{\Bb{R}}= (\frk{h}\cap\frk{k})\oplus(\frk{h}\cap\frk{p})$
$$\begin{array}{l} 
\rm{with}\qquad \qquad \qquad \quad \frk{h}\cap\frk{p} = \underset{\beta\in
\triangle^{re}}{\bigoplus}\Bb{R}(p_{_{\beta}}-p_{_{\sigma\beta}}), \qquad 
\qquad \qquad \qquad\\
{\rm{and}}\qquad \qquad \quad \qquad \,\, \frk{h}\cap\frk{k} = 
\big(\underset{\alpha\in\triangle^{re}}i\Bb{R}p_{\alpha}\big)\oplus 
\big(\underset{\beta\in\triangle^{re}}{\bigoplus}\Bb{R}(p_{_{\beta}} + 
p_{_{\sigma\beta}})\big) (\rm{modulo\ the\ center}), \qquad \qquad \qquad 
\qquad \end{array}$$ 
where $(p_\gamma)_{\gamma\in \triangle^{re}}\subset\frk{h}$ is the dual basis 
of the real roots. Clearly $(\frk{h}\cap\frk{k})\oplus i(\frk{h}\cap\frk{p})$ 
is a maximally compact Cartan subalgebra of $\frk{g}_{\Bb{R}}$. \end{rem}

\section{Vogan Diagrams}

\para \label{recall.1} Let $\frk{g}_{\Bb{R}}$ be an almost compact real form 
of $\frk{g}$ and $\sigma'$ be a semi-involution of second type  associated to 
$\frk{g}_{\Bb{R}}$. Let $\sigma$ be a Cartan involution adapted to $\sigma'$ 
and $\frk{g}_{\Bb{R}} = \frk{k}\oplus\frk{p}$ the corresponding Cartan 
decomposition of $\frk{g}_{\Bb{R}}$. Let 
$\frk{h}_{\Bb{R}} = (\frk{\dot{t}})^{\mu}\otimes1 \oplus 
(\frk{\dot{a}})^{\mu}\otimes1\oplus\Bb{R}ic\oplus\Bb{R}id$ be a maximally 
compact Cartan subalgebra of $\frk{g}_{\Bb{R}}$. By Lemma \ref{M.c}, 
$\frk{h}_{\Bb{R}}$ does not have any $\sigma$-real roots. Choose a positive 
system $\triangle^+$ for $\triangle(\frk{g},\frk{h})$, built from a basis of 
$i((\frk{\dot{t}})^{\mu}\otimes1)$ followed by a basis of 
$(\frk{\dot{a}})^{\mu}\otimes1$. Since 
$\sigma\vert_{i(\frk{\dot{t}})^{\mu}\otimes1} = Id$,
$\sigma\vert_{(\frk{\dot{a}})^{\mu}\otimes1} = -Id$ and $\frk{h}_{\Bb{R}}$
contains no $\sigma$-real roots, $\sigma(\triangle^+) =\triangle^+.$
Thus $\sigma$ fixes the $\sigma$-imaginary roots and permutes in 2-cycles
the $\sigma$-complex roots. 

A Vogan diagram of the triple 
($\frk{g}_{\Bb{R}},\frk{h}_{\Bb{R}},\triangle^+$), is a Dynkin diagram 
of $\triangle^+$ with the 2-element orbits of $\sigma$ so labeled and with the 
one element orbit painted or not, accordingly as the corresponding 
$\sigma$-imaginary simple root is noncompact or compact. In addition to this,
the underlying Dynkin has numerical labels as given in Figure 1(Section 1).  

\para\label{Int} The Cartan involution $\sigma$, is an involution of first 
kind. Let the realization $\frk{l}(\frk{\dot{s}},\zeta,\varepsilon_m)$ of 
$\frk{g}$ be adapted to $\sigma$. By paragraph {\bf{\ref{sigma}}}, $\sigma$ is 
of the form $\rho^{\jmath}\exp i\pi \ad(h_0)$, for 
$h_0\in\frk{h}^{\rho}_{\Bb{Z}}$, $\jmath \in\Bb{Z}_2$, where $\rho$ a diagram 
automorphism of $S(A)$. Then the base $\Pi = (\alpha_j)_{j\in[0,l]}$ can be 
chosen such that either $h_0=p_j\in\frk{h}$, for some $j\in[0,l]$ or 
$h_0=\ocirc{p_j}\in\ocirc{\frk{h}}\subset{\frk{h}}$, for some $j\in [1,l]$, 
where the set $(p_j)_{j\in[0,l]}\subset\frk{h}$ satisfies the property that 
$\alpha_k(p_j)=\delta_{k,j}$ for $k,j\in[0,l]$ (cf. \ref{p-i}) and 
$(\ocirc{p_j})_{j\in[1,l]}$ is a dual basis of the base 
$(\alpha_j)_{j\in[1,l]}$ of $\ocirc{\frk{g}}\subset{\frk{g}}$. 
Since $\alpha_k(p_j)=\delta_{k,j}$ for $k, j\in\lb0,1,\cdots,l\rb$, it can be 
easily seen that,  
\begin{equation}\label{eqp-i}
\begin{array}{ll}\exp i\pi\ad(p_k)e_j = (-1)^{\delta_{k,j}}e_k, \qquad& 
\text{for $j\in \lb 0,1,\cdots,l\rb$}.
\end{array}\end{equation} 
If the realization $\frk{l}(\frk{\dot{s}},\zeta,\varepsilon_m)$ 
of $\frk{g}$ is adapted to $\sigma = \exp i\pi\ad(p_0)$, then it follows from
the definition of adapted realization and discussions following 
Proposition \ref{H} that, 
$$\begin{array}{rlll}
\text{for $m=2$,} \quad & \dot{\frk{s}}=\dot{\frk{g}} & \zeta\otimes1(e_0) = 
-e_0 , & \quad \zeta\otimes1(e_j)=e_j,\ \text{for $j\neq 0$}, 
\qquad \qquad\qquad \qquad\\
\text{ for $m = 6$,}\quad & \dot{\frk{s}}=\dot{\frk{g}}\times\dot{\frk{g}} & 
\zeta^3\otimes1(e_0) = -e_0, & \quad \zeta\otimes1(e_j)=e_j,\  
\text{for $j\neq 0$.} \qquad \qquad\qquad \qquad \end{array}$$
By Lemma 3.12, $\frk{\dot{h}}^{\zeta}$ 
is the Cartan subalgebra of $\frk{l}(\frk{\dot{s}},\zeta,\varepsilon_m)''$
and $\exp i\pi \ad(p_0)\vert_{\frk{\dot{h}}^{\zeta}} =Id = 
(\zeta \otimes 1)\vert_{\frk{\dot{h}}^{\zeta}}$. Since $\lb e_j,f_j
\rb_{j\in[0,l]}$ and $\frk{\dot{h}}^{\zeta}$ generate $\frk{g}$, we get:
$$\begin{array}{ll} \exp i\pi\ad(p_o) = \zeta\otimes1, & \text{for}\ m=2,\quad 
\text{and} \quad \exp i\pi\ad(p_o) = \zeta^3\otimes1,  \text{for}\ m=6.
\end{array}$$
Also we have, $\alpha_k(\ocirc{p_j})=\delta_{k,j}$ 
for $k, j\in\lb1,\cdots,l\rb$  and $\alpha_0(\ocirc{p_j})=-a_j$ for 
$j\in \lb 1,\cdots,l\rb$. Since $p_0=d$ and $p_j=\ocirc{p_j}+a_jd$, for 
$j\in[1,l]$, using Equation (4.1) we get:
\begin{equation}\label{eq0p-i}\begin{array}{ll}
\exp i\pi\ad(\ocirc{p_k})e_j = (-1)^{\delta_{k,j}} e_k, \qquad& \text{for 
$k, j\in \lb 1,\cdots,l\rb$,}\\ 
\exp i\pi\ad(\ocirc{p_k})e_0 = (-1)^{a_k}e_0, 
\qquad& \text{for $k\in\lb1,\cdots,l\rb$},\\
\end{array}\end{equation}
Thus it follows from Eqn (\ref{eqp-i}) and Eqn (\ref{eq0p-i}) that for
$k\in \lb 1,2,\cdots,l\rb$, 
\begin{equation}
\begin{array}{ll}\exp i\pi {\rm{ad}}(p_k) = \exp i \pi \ad(\ocirc{p_k}),  &
\text{whenever $a_k$ is even for $k\in\lb1,\cdots,l\rb$ },\\
\exp i\pi\ad(p_k) = (\exp i \pi \ad(\ocirc{p_k})\zeta)\otimes1,  &
\text{whenever $a_k$ is odd for $k\in\lb1,\cdots,l\rb$ }.\end{array} 
\end{equation}

From the symmetry of the diagram it is clear that for $\frk{g} = 
\frk{l}(\frk{\dot{g}}, \mu,\varepsilon_k)$, the involution of first type 
$\sigma$, is of the form $\rho \exp i\pi\ad(p_j)$, $j\in [1,l]$, only when 
$k=2$. In this case $\frk{l}(\frk{\dot{g}}, \zeta, \varepsilon_{2k})$ =
$\frk{l}(\frk{\dot{g}}, \mu \exp \frac{i\pi}{2}\ad(\ocirc{p_{\rho_{(0)}}}),
\varepsilon_{2k})$ is an adapted realization of 
$\frk{g}$ where $\ocirc{p_{\rho_{(0)}}}$ is the dual of $\rho(\alpha_0)$ in 
$(\frk{\dot{h}})^{\mu}$ and  one gets that,
\begin{equation}
\begin{array}{ll}\exp i\pi{\rm{ad}}(p_j) = \exp i \pi \ad(\ocirc{p_j}),  &
\text{whenever $a_j$ is even for $j\in\lb1,\cdots,l\rb$ },\\
\exp i\pi\ad(p_j) = (\exp i \pi \ad(\ocirc{p_j})\zeta^2)\otimes1,  &
\text{whenever $a_j$ is odd for $j\in\lb1,\cdots,l\rb$ }.\end{array} 
\end{equation}

\para \label{Cn} It follows from Eqn(\ref{eqp-i}) that the \aV diagram 
associated to the involution $\exp i\pi\ad(p_j)$, $j\in[0,l]$ has exactly one 
painted vertex, namely the $i^{th}$ vertex. By  Eqn(\ref{eq0p-i}), the \aV 
diagrams associated to the involution $\exp i\pi\ad(\ocirc{p_j})$, $j\in[1,l]$ 
have exactly one painted vertex, namely the $j^{th}$ vertex, if $a_j$ is even 
and exactly two painted vertices, namely the $j^{th}$ and the $0^{th}$ 
vertices, if $a_j$ is odd.  Note that $\sigma$ is of the form $\rho H$ for 
$\rho\neq Id$, only when $\frk{g}$ is of type $A_{2l-1}^{(2)}$ or 
$D_{l+1}^{(2)}$. In both the cases $\alpha_0$ is a $\sigma$-complex root and 
hence the interior automorphism $H$ is of the form $\exp i\pi\ad(p_j)$, for 
some $j\in[1,l]$ such that $\rho(\alpha_j) =\alpha_j $.

\para {\it{Equivalence of \aV diagrams}}: 
In \cite[Chapter VI, Ex.18]{Knapp}, an operation $R[j]$ on the Vogan diagram of
a simple finite dimensional Lie algebra $\ocirc{\frk{g}}$ is defined as 
follows: \\$R[j]$ acts on the base $\ocirc{\pi}$ of $\ocirc{\frk{g}}$ by 
reflection corresponding to the noncompact simple root $\alpha_j$. As a 
consequence of $R[j]$, the colour of $\alpha_j$ and all vertices not adjacent 
to $\alpha_j$ remain unchanged; if $\alpha_k$ is joined to $\alpha_j$ by a 
double edge and $\alpha_j$ is the smaller root, then the colour of $\alpha_k$ 
remains unchanged and if $\alpha_k$ is joined to $\alpha_j$ by single or triple
lines then the colour of $\alpha_k$ is reversed. 
 
Let the following operations generate equivalence relation between \aV 
diagrams:\\ 
1.\ Application of a diagram automorphism of the Dynkin diagram $S(A)$.\\ 
2.\ Application of a sequence of $R[j]$s for $j\in\lb1,2,\cdots,l\rb$  if
the $j^{th}$ vertex is coloured. 

Given a Lie algebra $\frk{g} = \frk{l}(\frk{\dot{g}},\mu,\varepsilon_k)$, the 
set of involutions $\lb\exp i\pi \ad(\ocirc{p_k})\rb_{k\in[1,l]}$ is a subset 
of $Aut((\frk{\dot{g}})^{\mu})$. Observe that the base of the twisted affine 
Kac-Moody Lie algebra $\frk{g}$ is given by $\ocirc{\pi}\cup\lb\alpha_0\rb$, 
where $\alpha_0=\delta-\theta$, for $\frk{g}$ not of type $A_{2l}^{(2)}$, and
$\alpha_0=\delta-2\theta$, for $\frk{g}$ of type $A_{2l}^{(2)}$, 
$\theta \in \ocirc{\triangle}_+$;  $\ocirc{\pi}$ is a base of 
$(\frk{\dot{g}})^{\mu}$ and $\ocirc{W}\ltimes T$ is the the Weyl group of 
$\frk{g}$ (cf.Section1) with $w(\delta)=\delta$ for $w\in W$. Hence, if the 
compact real forms of $\ocirc{\frk{g}}$ corresponding to the involutions 
$\exp i\pi\ad(\ocirc{p_k})$ and $\exp i\pi\ad(\ocirc{p_j})$ are isomorphic via 
an isomorphism induced by $\ocirc{W}$, then $\ocirc{W}\subset W$ induces an 
isomorphism between the corresponding almost compact real forms of $\frk{g}$. 
As a consequence we get the second equivalence relation as defined above.

The following are some examples of equivalent \aV diagrams. \\
Example 1. $D_{l+1}^{\left(2\right)} \left(l \geq 2 \right)$
$$\begin{array}{lcl}
\begin{array}{l}\underset{\alpha_{0}}{\stackrel{\overset{\qquad}{1}}{\bullet}} 
\Leftarrow  \underset{\alpha_{1}} {\stackrel{1}{\circ}}  - \cdots- 
  \underset{\alpha_{l-1}}{\stackrel{1}{\circ}}  \Rightarrow
\underset{\underset{\qquad}{\alpha_{l}}}{\stackrel{1}{\circ}}\\
\quad \sigma = \exp i\pi\ad(p_0)
\end{array} &\stackrel{\text{equivalence relation 1 }}{\leftrightarrow}&
\begin{array}{l}\underset{\alpha_{0}}{\stackrel{\overset{\qquad}{1}}{\circ}} 
\Leftarrow  \underset{\alpha_{1}} {\stackrel{1}{\circ}}  - \cdots- 
  \underset{\alpha_{l-1}}{\stackrel{1}{\circ}}  \Rightarrow
\underset{\underset{\qquad}{\alpha_{l}}}{\stackrel{1}{\bullet}}\\
\quad \sigma = \exp i\pi\ad(p_l)
 \end{array} \end{array}$$ 
Example 2. $A_{2l-1}^{(2)}$
$$\begin{array}{lcl}
\begin{array}{l} \, \qquad \stackrel{1}{\bullet} _{\alpha_0}\\
\, \qquad \, \vert \\
   \underset{\alpha_{1}}{\stackrel{1}{\bullet}} -
     \underset{\alpha{2}}{\circ^{2}} -\cdots
   \underset{\alpha_{i}}{\stackrel{2}{\circ}} - \cdots-
   \underset{\alpha_{l-1}}{\stackrel{2}{\circ}} \Leftarrow
   \underset{\alpha_{l}}{\stackrel{1}{\circ}}\\
\quad \sigma = \exp i\pi\ad(\ocirc{p_1})
 \end{array}&  \stackrel{R[l-1]\cdots R[2]R[1]}
{\underset{\text{equivalence relation 2}}{\leftrightarrow}}  &
\begin{array}{l} \, \qquad \stackrel{1}{\circ} _{\alpha_0}\\
\, \qquad \, \vert \\
   \underset{\alpha_{1}}{\stackrel{1}{\circ}} -
     \underset{\alpha{2}}{\circ^{2}} -\cdots
   \underset{\alpha_{i}}{\stackrel{2}{\circ}} - \cdots-
   \underset{\alpha_{l-1}}{\stackrel{2}{\bullet}} \Leftarrow
   \underset{\alpha_{l}}{\stackrel{1}{\circ}}\\
\quad \sigma = \exp i\pi\ad(\ocirc{p_{l-1}})
 \end{array}\end{array}$$ 
Example 3. $A_{2l-1}^{(2)}$
$$\begin{array}{lcl}
\begin{array}{l}
\overset{\alpha_0}{_1\circ}  \\
\overset{\diagdown}
{\underset{\diagup} {\updownarrow \quad \circ_{\alpha_{_2}}^2}}
- \underset{\alpha_{3}}{\stackrel{2}{\bullet}} - 
- \underset{\alpha_{4}}{\stackrel{2}{\circ}} - 
   \underset{\alpha_{5}}{\stackrel{2}{\circ}} \Leftarrow
   \underset{\alpha_{6}}{\stackrel{1}{\circ}}\\
\underset{\alpha_1}{_1\circ}\end{array} &
 \stackrel{R[4]R[3]R[4]R[2]R[4]R[3]}{\underset{\text{Eq.rel 2}}{\leftrightarrow}}  &
\begin{array}{l}
\overset{\alpha_0}{_1\circ}  \\
\overset{\diagdown}
{\underset{\diagup} {\updownarrow \quad \circ_{\alpha_{_2}}^2}}
- \underset{\alpha_{3}}{\stackrel{2}{\circ}} - 
- \underset{\alpha_{4}}{\stackrel{2}{\bullet}} - 
   \underset{\alpha_{5}}{\stackrel{2}{\circ}} \Leftarrow
   \underset{\alpha_{6}}{\stackrel{1}{\circ}}\\
\underset{\alpha_1}{_1\circ}\end{array}
\end{array}$$ 
Example 4.$A_{2l-1}^{(2)}$
$$\begin{array}{lclcl}
\begin{array}{l} \, \qquad \stackrel{1}{\bullet} _{\alpha_0}\\
\, \qquad \, \vert \\
   \underset{\alpha_{1}}{\stackrel{1}{\circ}} -
     \underset{\alpha_{2}}{\bullet^{2}} -\cdots
   \underset{\alpha_{l-1}}{\stackrel{2}{\circ}} \Leftarrow
   \underset{\alpha_{l}}{\stackrel{1}{\circ}}
 \end{array}&  \stackrel{\text{Eq.rel 1}}{\leftrightarrow}  &
\begin{array}{l} \, \qquad \stackrel{1}{\circ} _{\alpha_0}\\
\, \qquad \, \vert \\
   \underset{\alpha_{1}}{\stackrel{1}{\bullet}} -
     \underset{\alpha_{2}}{\bullet^{2}} -\cdots
   \underset{\alpha_{l-1}}{\stackrel{2}{\circ}} \Leftarrow
   \underset{\alpha_{l}}{\stackrel{1}{\circ}}\\
 \end{array}& \stackrel{R[1]}{\underset{\text{Eq.rel 2}}{\leftrightarrow}}  &
\begin{array}{l} \, \qquad \stackrel{1}{\circ} _{\alpha_0}\\
\, \qquad \, \vert \\
   \underset{\alpha_{1}}{\stackrel{1}{\bullet}} -
     \underset{\alpha_{2}}{\circ^{2}} -\cdots
   \underset{\alpha_{l-1}}{\stackrel{2}{\circ}} \Leftarrow
   \underset{\alpha_{l}}{\stackrel{1}{\circ}}\\
 \end{array}
\end{array}$$ 
In Example 4, ``the equivalence relation 2'' gives $r_1(\alpha_2) = 
\alpha_2+\alpha_1$ and $r_1(\alpha_0) = \alpha_0$ and hence the final diagram 
is a consequence of the fact that $[\frk{k},\frk{k}]\subset\frk{k}, 
[\frk{k},\frk{p}]\subset\frk{p}$ and $[\frk{p},\frk{p}]\subset\frk{k}$.

\begin{rem} \label{painted-roots}As a consequence of paragraph {\bf{\ref{Int}}}
and the equivalence relation between the \aV diagrams, there exists a base 
$\Pi'$ of a twisted affine Kac-Moody Lie algebra $\frk{g}$ with corresponding 
positive root system $\triangle'$ such that the \aV diagram associated to 
$(\frk{g}_{\Bb{R}},\frk{h}_{\Bb{R}},\triangle')$ has at most two painted simple
roots. Hence every \aV diagram is equivalent to one with at most two painted
simple roots. \end{rem}

\begin{defn} An abstract \aV diagram is an irreducible abstract Dynkin 
diagram of a twisted affine Kac-Moody Lie algebra $\frk{g}$ indicated
with the following additional structures:\\
1. A diagram automorphism of order 1 or 2, indicated by labeling the 2-element 
orbits. \\
2. A subset of the 1-element orbits, indicated by painting the vertices 
corresponding to the members of the subset.\\
3. The vertices of the Dynkin diagram have enumeration as given in Figure 1.
\end{defn}
\noindent Every \aV diagram is of course an abstract \aV diagram. 

\begin{thm}\label{abs}
Let an abstract \aV diagram for a twisted affine Kac-Moody Lie 
algebra be given. Then there exists an almost compact real form of a twisted 
affine Kac-Moody Lie algebra such that the given diagram is the \aV diagram
of this almost compact real form. \end{thm}
\pf It is known from \cite{Kac}, that a GCM $A$ can be uniquely associated to 
Dynkin diagram $S(A)$ and its enumeration of vertices. By 
\cite[Proposition 1.1]{Kac}, there exists a unique up-to-isomorphism 
realization for every GCM $A$. Following \cite[Sect. 1.3]{Kac}, one can 
associate with A, a Kac-Moody Lie algebra  $\frk{g}$ = $\frk{g}(A)$.

By Proposition \ref{1-2pair}, there exists a one to one correspondence 
between involutions of first kind $\sigma$ and the pairs ($\sigma',\frk{h}$) 
formed of a semi-involution of second type $\sigma'$ and a maximally compact 
Cartan subalgebra for $\sigma'$. Thus if the involution of 
$\frk{g}(A)$ of the first kind $\sigma$, can be extracted from the additional 
structural information superimposed on the \aV diagram, then one can 
associate to the \aV diagram an almost compact real form of a twisted 
affine Kac-Moody Lie algebra. We shall now prescribe an algorithm to associate 
an involution of the first kind $\sigma$ to a given \aV diagram.

Let $V(A)$ denote a \aV diagram of $\frk{g}(A)$. By Remark 
\ref{painted-roots}, $V(A)$ is equivalent to a  \aV diagram of $\frk{g}(A)$
with at most two painted simple roots.\\
1. If $V(A)$ has no painted simple roots and no 2-element orbits, then it 
corresponds the compact form $\frk{u}(A)$ of $\frk{g}(A)$ and $\sigma =Id$ in 
this case. \\
2. Suppose $V(A)$ contains no 2-element orbit. Then:\\
$\bullet$ $\sigma =\exp i\pi\ad(p_j)$, for $j\in[0,l]$, if  only the $j^{th}$ 
vertex is painted.\\
$\bullet$ $\sigma =\exp i\pi\ad(\ocirc{p_j})$, for $j\in[1,l]$, if $a_j$ is odd
and the $0^{th}$ and $j^{th}$ vertices are painted.\\
3. Suppose $V(A)$ contains 2-element orbits. If $\rho$ denotes the Dynkin 
diagram automorphism in $V(A)$ then from $\S 4.3$ it follows that:\\
$\bullet$ $\sigma =\rho\exp i\pi\ad(p_i)$, for $i\in[1,l]$, if  only the 
$i^{th}$ vertex is painted.\\
$\bullet$ $\sigma =\rho$, if no vertices are painted.\\The association of 
$\sigma$ with $V(A)$ thus completes the proof of the theorem.
$\hfill \blacksquare$ \endpf

\begin{thm}\label{eq} 
Two almost compact real forms of a twisted affine Kac-Moody Lie 
algebra $\frk{g}$ having equivalent \aV diagram are isomorphic. \end{thm}
\begin{Rem} Owing to the definition of the equivalence relation 
between \aV diagrams, to prove the theorem, it suffices to show that two 
almost compact real forms having the same \aV diagram are isomorphic.
\end{Rem}
\pf
\noindent Let $\frk{g}_1$ and $\frk{g}_2$ be two almost compact real forms of 
$\frk{g}$ having the same \aV diagram. As they both have the same Dynkin 
diagram with the same enumeration on the vertices, the same generalized Cartan 
matrix $A$ is associated with both $\frk{g}_1$ and $\frk{g}_2$. Thus the 
unique twisted Kac-Moody Lie algebra $\frk{g}=\frk{g}(A)$ is the 
complexification of $\frk{g}_1$ and $\frk{g}_2$. By Proposition \ref{CSI}
there exists Cartan semi-involutions $\vartheta_1, \vartheta_2$ adapted to
$\frk{g}_1$ and $\frk{g}_1$ respectively. Let $\frk{u}_j = 
\frk{g}^{\vartheta_j}$, for $j=1,2$, be the corresponding compact real forms 
of $\frk{g}$. Let $\sigma$ be the involution represented by the \aV 
diagram. Then for $j=1,2$, $\sigma\vert_{\frk{u}_j} = \varpi_j$ is the 
corresponding Cartan involution on $\frk{u}_j$. Since by \cite[Theorem 4.6]{R4}
all Cartan semi-involutions are conjugate by Int$(\frk{g})$, there exits 
$x\in$ Int($\frk{g}$) such that $x\vartheta_1x^{-1}=\vartheta_2$ and 
consequently $x.\frk{u}_1=\frk{u}_2$. As $x.\frk{g}_1$ is isomorphic to 
$\frk{g}_1$, without loss of generality we may assume from the outset that
$\frk{u}_1 = \frk{u}_2 = \frk{u}$ and we have  
$$ \qquad \qquad\qquad    \varpi_j(\frk{u}) = \frk{u},
\quad \quad\quad \quad \rm{for} \quad j=1,2 .$$
Let $\frk{h}_1 = \frk{t}_1 \oplus \frk{a}_1$ and $\frk{h}_2 = \frk{t}_2 \oplus 
\frk{a}_2$ be the Cartan decompositions of the Cartan subalgebras of $\frk{g}_1
$ and $\frk{g}_2$ respectively, where  $\frk{t}_j$ and $\frk{a}_j$, for 
$j=1,2$, are respectively the +1 and -1 eigenspaces of $\sigma$ in $\frk{h}_j$.
Consequently, for $j=1,2$, $\frk{t}_j\oplus i\frk{a}_j$ is a maximal abelian 
subspace of $\frk{u}$ and by Remark \ref{h.decomp}, 
$\frk{t}_j\oplus i\frk{a}_j$, is a  maximally compact Cartan subalgebra of 
$\frk{u}$. Hence by \cite[Proposition 4.9c]{R4},$\frk{t}_1\oplus i\frk{a}_1$ 
and $\frk{t}_2\oplus i\frk{a}_2$ are conjugate by an element $k\in$ 
Int$(\frk{u})$. Replacing $\frk{g}_2$ by $k\frk{g}_2$ and arguing as above we 
may assume that $\frk{t}_1\oplus i\frk{a}_1 =\frk{t}_2\oplus i\frk{a}_2$. Thus 
$\frk{t}_1\oplus i\frk{a}_1$ and $\frk{t}_2\oplus i\frk{a}_2$ have the same 
complexification, which is denoted by $\frk{h}$.

Now the complexifications $\frk{g}$ and $\frk{h}$ have been aligned and the 
root systems are the same. Let the respective positive root systems be given by
$\triangle_+^1$ and $\triangle_+^2$. By paragraph {\bf{\ref{W(u)}}} there 
exists an interior automorphism $s \in$ Int$(\frk{u})$ such that 
$s \triangle_+^2= \triangle_+^1$. Again replacing $\frk{g}_2$ by 
$s.\frk{g}_2$ and repeating above argument we assume 
$\triangle_+^1 = \triangle_+^2 = \triangle_+$ from the outset. 

Using the conjugacy of the compact real forms of 
$\frk{g}$ we construct in this case,
$$\begin{array}{l} \frk{u} = 
\underset{\alpha \in \ocirc{\triangle}}{\sum} 
\Bb{R}iH_{\alpha} \oplus \underset{\gamma \in \triangle^{re}}{\sum}
\Bb{R}i\big(e_{\gamma}+f_{\gamma}\big)
\oplus \underset{\gamma \in \triangle^{re}}{\sum} 
\Bb{R}\big(e_{\gamma}-f_{\gamma}\big) \oplus  \Bb{R}ic \oplus \Bb{R}id \\ 
\qquad \qquad \oplus \underset{\alpha \in \ocirc{\triangle}_l}{\sum}
\Bb{R}i\big(H_{\alpha}t^{kn}+H_{\alpha}t^{-kn}\big) \oplus 
\underset{\alpha \in \ocirc{\triangle}_l}{\sum} \Bb{R}\big(
H_{\alpha}t^{kn}-H_{\alpha}t^{-kn}\big)\oplus \\
\qquad \big(\bigoplus_{j=0}^{k-1} \Big( \sum_{\alpha \in \ocirc{\triangle}_s} 
\Bb{R}i\big(h_{\alpha}^{^j}t^{kn+j}+h_{\alpha}^{^j}t^{-kn-j}\big) \oplus 
\sum_{\alpha \in \ocirc{\triangle}_s} \Bb{R}\big(h_{\alpha}^{^j}t^{kn+j} - 
h_{\alpha}^{^j}t^{-kn-j}\big)\Big)\big),\end{array}$$ where 
for $\gamma\in\triangle^{re}$, 
$e_{\gamma} = E_{\alpha,ks}$, $f_{\gamma}=E_{-\alpha,-ks}$ for $\alpha\in
\ocirc{\triangle}_l$, $s\in\Bb{Z}$ and $e_{\gamma}$ = $E_{\alpha,ks+j}$, 
$f_{\gamma}= E_{-\alpha,-ks-j}$ for $\alpha\in\ocirc{\triangle}_s$, 
$s\in\Bb{Z}$ and $j\in \Bb{Z}_k$. Here, $\lb E_{\alpha,ks},E_{-\alpha,-ks},
H_{\alpha},\ \text{for}\ \alpha\in\ocirc{\triangle}_l \rb$ and 
$\lb E_{\alpha,ks+j}, E_{-\alpha,-ks-j}, h_{\alpha}^{^j},\ \text{for}\ 
\alpha\in\ocirc{\triangle}_s, j\in \Bb{Z}_k\rb $ are defined as in 
paragraph {\bf{\ref{long/short}}}. 

Case 1: Suppose $\frk{h}_{\Bb{R}}$ has no $\sigma$-complex roots.
As the \aV diagram for $\frk{g}_1$ and $\frk{g}_2$ are the same, the 
automorphisms of $\triangle_+$ defined by $\varpi_{1}$ and $\varpi_{2}$ 
have the same effect on $\frk{h}^{*}$. Thus, 
\begin{eqnarray}  
\varpi_{1}(H) = \varpi_{2}(H),\qquad \text{for all}\ H\in 
\frk{h},\\ \varpi_1(\triangle^{im}) = \varpi_2(\triangle^{im}), \quad \qquad
\quad\qquad\   \\ \varpi_1\vert_{\Bb{R}ic\oplus\Bb{R}id} = 
\varpi_2\vert_{\Bb{R}ic\oplus\Bb{R}id}.\qquad\quad\quad\quad \ 
\end{eqnarray} 
If $\alpha$ is a simple $\sigma$-imaginary real
root then \begin{eqnarray}   \varpi_1(e_j) = e_j =\varpi_2(e_j), &
\quad \ \text{if the}\ j^{th} \ \text{vertex is unpainted},\\
\varpi_1(e_j) = - e_j =\varpi_2(e_j), & \text{if the}\ j^{th}\ \text{vertex
is painted}.\end{eqnarray} 
Since $\frk{h}$ and $\lb e_j,f_j, j\in [0,l]\rb$  generate $\frk{g}$, it 
follows that $\varpi_1 = \varpi_2$ on $\frk{u}$, hence $\frk{k}_1 = 
\frk{u}^{\varpi_1} = \frk{k}_2$ and if for $j=1,2$, $\frk{p}_j$ denotes the -1
eigenspace of $\varpi_j$ on $\frk{u}$, then $\frk{p}_1 =\frk{p}_2$. 
Hence $$\frk{g}_1 = 
\frk{k}_1\oplus i\frk{p}_1 = \frk{k}_2\oplus i\frk{p}_2=\frk{g}_2.$$

Case 2: Suppose there exists $\sigma$-complex simple real roots in $\frk{g}_1$
and $\frk{g}_2$. Let $\rho$ denote the diagram automorphism of $S(A)$. In this 
case the $\sigma$-imaginary roots are treated as in Case 1. If for all 
$\sigma$-complex roots $\gamma\in\triangle^{re}$ there exists 
$X_{\gamma}\in\frk{g}_{\gamma}$,$X_{\rho\gamma}\in\frk{g}_{\rho\gamma}$
and constants $a_{\gamma},b_{\gamma}$ such that
\begin{equation}\label{w12}
\varpi_1(X_{\gamma}) = a_{\gamma}\ X_{\rho \gamma} \qquad \text{and}
\qquad \varpi_2(X_{\gamma}) = b_{\gamma}\ X_{\rho \gamma},
\end{equation} then the same calculations as done in 
\cite[Theorem 5.2(Case 2)]{B1} show that for each pair of $\sigma$-complex
simple roots $\gamma, \varpi\gamma$, square roots 
$a_{\gamma}^{1/2}, a_{\rho\gamma}^{1/2},$ $b_{\gamma}^{1/2}, 
b_{\rho\gamma}^{1/2}$ can be chosen such that 
$a_{\gamma}^{1/2}, a_{\rho\gamma}^{1/2} =1,$ and 
$b_{\gamma}^{1/2}, b_{\rho\gamma}^{1/2}=1.$ Then defining
$H, H'\in \frk{u}\cap\frk{h}$ such that $\alpha(H)=\alpha(H')=0$ for $\alpha$
$\sigma$-imaginary simple root and 
$$\begin{array}{ll}
\exp (\frac{1}{2}\gamma(H))=a_\gamma^{1/2}, &
\exp (\frac{1}{2}\rho\gamma(H))=a_{\rho\gamma}^{1/2},\\
\exp (\frac{1}{2}\gamma(H'))=b_\gamma^{1/2}, &
\exp (\frac{1}{2}\rho\gamma(H'))=b_{\rho\gamma}^{1/2},
\end{array}$$
for $\gamma, \ \rho\gamma$ $\sigma$-complex simple roots, 
similar calculations as in \cite[Theorem 5.2(Case 2)]{B1} show that application
of the identity $$ \varpi_2 \circ {\rm{Ad}}(\exp\frac{1}{2}(H-H')) =
{\rm{Ad}}(\exp\frac{1}{2}(H-H')) \circ \varpi_1,$$ to the $\sigma$-eigenspaces
of $\frk{g}_1$ gives an isomorphism between $\frk{g}_1$ and $\frk{g}_2$.

Thus to complete the proof of the theorem we need to show the existence of
$X_{\gamma}\in\frk{g}_{\gamma}$, $X_{\rho\gamma}\in\frk{g}_{\rho\gamma}$
and constants $a_{\gamma},b_{\gamma}$ 
satisfying Eqn(\ref{w12}). Observe from Figure 1, that for $\frk{g}$ of type 
Aff $k, k\neq 1$, $\sigma$-complex roots exist only when $\frk{g}$ of type 
$A_{2\ell-1}^{(2)}$ or $D_{\ell+1}^{(2)}$. 
If $\gamma=\alpha+ks\delta\in\triangle^{re}$ for 
$\alpha\in \ocirc{\triangle}_l$, then $\frk{g}_{\gamma}=\Bb{C}e_{\alpha}\otimes
t^{ks}$, where $e_{\alpha}\in (\dot{\frk{g}})_{\alpha}$.
As $e_{\alpha}\otimes t^{ks}\in \frk{l}(\frk{\dot{g}},Id,1)$, by 
\cite[Theorem 5.2(Case 2)]{B1}, Eqn (\ref{w12}) is satisfied in this case.
However problems can arise if a $\sigma$-complex root $\gamma\in\triangle^{re}$
is a short root. Note that the following are the only short $\sigma$-complex 
simple roots:
\begin{equation}\label{rho}
\begin{array}{ll} \alpha_0,\ \alpha_1=\alpha_{\rho0}, \qquad &\qquad  
\text{when $\frk{g}$ of type $A_{2\ell-1}^{(2)}$},\\ 
\alpha_0,\ \alpha_\ell=\alpha_{\rho0}, \qquad &\qquad  \text{when $\frk{g}$ of
type $D_{\ell+1}^{(2)}$}. \end{array}\end{equation} 
Let $H_{\rho_{(0)}}$:= $\exp \frac{i\pi}{2}\ad(\ocirc(p_{\rho_{(0)}}))$. Then
$\frk{l}(\frk{\dot{g}}, \mu H_{\rho_{(0)}}, \varepsilon_4)$ is a realization
of $\frk{g}$ adapted to the involution $\sigma$, having a 2-element orbit.
If $\frk{g}_{\alpha_0} = \Bb{C}X_{\alpha_0}$ and $\frk{g}_{\rho(\alpha_0)} = 
\Bb{C}X_{\alpha_j}$, for $\alpha_j$ a simple short root of $\frk{g}$, then in 
the realization, $\frk{l}(\frk{\dot{g}}, \mu H_{\rho_{(0)}}, \varepsilon_4)$
we have,
\begin{equation}\label{base}\begin{array}{ll}
X_{\alpha_{_0}} = (e_{-\theta^0}-e_{-\mu\theta^0})\otimes t & \text{for}\ 
\frk{g}\ \text{of type}\ D_{\ell+1}^{(2)}\ \rm{and}\ A_{2\ell-1}^{(2)}\\
X_{\alpha_j} =(e_{\dot{\alpha}_j}+e_{-\mu\dot{\alpha}_j})\otimes t \qquad 
&\text{where} \left\{ \begin{array}{l} j=\ell  \ \rm{for}\ \frk{g}\  
\text{of type}\ D_{\ell+1}^{(2)},\ \dot{\alpha_j}\in\triangle(D_{\ell+1})\\ 
j=1\ \rm{for}\ \frk{g}\ \text{of type}\ 
A_{2\ell-1}^{(2)}, \ \dot{\alpha_j}\in\triangle(A_{2\ell-1}) \end{array}\right.
\end{array}\end{equation}

Since for appropriate $j$ (as explained in Eqn(\ref{base})),
$\varpi_1(\frk{g}_{\alpha_0})\subset\frk{g}_{\alpha_j}$ and
$\varpi_2(\frk{g}_{\alpha_0})\subset\frk{g}_{\alpha_j}$, by Eqn(\ref{base})
there exists constants $a^1_{\theta^{0}}, a^2_{\theta^{0}}$, 
$a^1_{\mu\theta^{0}},a^2_{\mu\theta^{0}}$ and $b^1_{\theta^{0}}, 
b^2_{\theta^{0}}$, $b^1_{\mu\theta^{0}},b^2_{\mu\theta^{0}}$  such that
$$\begin{array}{ll} 
\varpi_1(e_{-\theta^0}\otimes t) = a^1_{\theta^{0}}  e_{\alpha_j}\otimes t +
a^1_{\mu\theta^{0}}\ e_{\mu\alpha_j}\otimes t, 
 \quad & \quad \varpi_1(e_{-\mu\theta^0}\otimes t) = 
a^2_{\theta^{0}}\ e_{\alpha_j}\otimes t + a^2_{\mu\theta^{0}} \ 
e_{\mu\alpha_j}\otimes t,\\ 
\varpi_2(e_{-\theta^0}\otimes t) = b^1_{\theta^{0}}  e_{\alpha_j}\otimes t +
b^1_{\mu\theta^{0}}\ e_{\mu\alpha_j}\otimes t, 
 \quad & \quad \varpi_2(e_{-\mu\theta^0}\otimes t) = 
b^2_{\theta^{0}}\ e_{\alpha_j}\otimes t + b^2_{\mu\theta^{0}} \ 
e_{\mu\alpha_j}\otimes t.
\end{array}$$
Claim: $a^1_{\theta^{0}}-a^2_{\theta^{0}} = 
a^1_{\mu\theta^{0}}-a^2_{\mu\theta^{0}}$; \, $b^1_{\theta^{0}}-b^2_{\theta^{0}}
= b^1_{\mu\theta^{0}}-b^2_{\mu\theta^{0}}$.

\noindent{\it{Proof the claim}}: Recall that $\varpi_1$ and $\varpi_2$ are
restrictions of the involution $\sigma$ to the compact forms $\frk{u}_1$ and
$\frk{u}_2$ adapted to $\frk{g}_{1}$ and $\frk{g}_{2}$ respectively.
Since for $j=1,2$, the Cartan involution $\sigma$ is adapted to both 
$\frk{g}_j$ and the realization 
$\frk{l}(\frk{\dot{g}},\mu H_{\rho_{(0)}}, \varepsilon_4)$ of $\frk{g}$ is 
adapted to ($\frk{g}_{j},\sigma,\frk{h}$), by \cite[Proposition 3.5]{M-R} there
exists $\dot{\sigma}_j\in$ Aut($\frk{\dot{g}}$) commuting with 
$\mu H_{\rho_{(0)}}$ such that $\varpi_j = \dot{\sigma}_j\otimes1$ on the 
compact form $\frk{u}$ of 
$\frk{l}(\frk{\dot{g}},\mu H_{\rho_{(0)}},\varepsilon_4)$.

As $\mu H_{\rho_{(0)}}(e_{\mu \theta^0}\otimes t) = -i\ e_{\theta^0}\otimes t$,
we have $\ e_{\theta^0}\otimes t = i\ 
\mu H_{\rho_{(0)}}(e_{\mu \theta^0}\otimes t)$. Hence,
$$\begin{array}{ll} 
a^2_{\theta^{0}}\ e_{\alpha_j}\otimes t + a^2_{\mu\theta^{0}} \ 
e_{\mu\alpha_j}\otimes t &= \varpi_1(e_{-\mu\theta^0}\otimes t)
=\varpi_1\big(i \mu H_{\rho_{(0)}} (e_{-\theta^0})\otimes t\big)\\
& = \dot{\sigma}_1\big(i(\mu H_{\rho_{(0)}}(e_{-\theta^0})\big)\otimes t
 = i\  \mu H_{\rho_{(0)}}\big(\dot{\sigma}_1(e_{-\theta^0})\big)\otimes t\\
&= i\  \mu H_{\rho_{(0)}} (\big (\dot{\sigma}_1(e_{-\theta^0}))\otimes t \big)
= i\  \mu H_{\rho_{(0)}} (\big(\varpi_1(e_{-\theta^0}\otimes t)\big)\\ 
&= i\  \mu H_{\rho_{(0)}} \big(a^1_{\theta^{0}} e_{\alpha_j}\otimes t +
a^1_{\mu\theta^{0}}\ e_{\mu\alpha_j}\otimes t \big) \\
&= i^2\  (a^1_{\theta^{0}}  e_{\mu\alpha_j}\otimes t + 
a^1_{\mu\theta^{0}}\ e_{\alpha_j}\otimes t) = 
-a^1_{\theta^{0}}  e_{\mu\alpha_j}\otimes t -
a^1_{\mu\theta^{0}}\ e_{\alpha_j}\otimes t.
\end{array}$$ Comparing coefficients we get, $a^2_{\theta^{0}} =
-a^1_{\mu\theta^{0}}$ and $a^2_{\mu\theta^{0}}= -a^1_{\theta^{0}}$. Hence 
$a^1_{\theta^{0}}-a^2_{\theta^{0}}= -a^2_{\mu\theta^{0}} +a^1_{\mu\theta^{0}}$
as desired. It can be similarly shown that $b^1_{\theta^{0}}-b^2_{\theta^{0}}= 
-b^2_{\mu\theta^{0}} +b^1_{\mu\theta^{0}}$. Thus for short $\sigma$-complex 
simple roots $X_{\alpha_0}\in\frk{g}_{\alpha_0}$ and 
$X_{\alpha_j}\in\frk{g}_{\rho\alpha_0}$ exists such that Eqn(4.10) is 
satisfied. Hence the theorem. $\hfill 
\blacksquare$ \endpf

Using the definition of equivalence relations between the \aV diagrams 
(cf. {\bf{4.4}}) and \cite[Figure 6.1,Figure 6.2]{Knapp} we give in the 
following table the non-equivalent \aV diagrams of the twisted affine Kac-Moody
Lie algebras corresponding to non-trivial involutions of first type. Note that 
owing to the equivalence relation as described in Example 1, the 
non-equivalent \aV diagrams for $\frk{g}$ correspond to the 
following involutions of first type : $\mu\otimes 1$, $\exp\ i\pi\ad(p_i)$ for 
$1\leq i\leq \frac{l}{2}$, $\exp\ i\pi\ad(\ocirc{p_i})$ for $1\leq i\leq l$.
Likewise, owing to the equivalence relations described in Examples 2 and 4, 
the non-equivalent \aV diagrams for $\frk{g}$ correspond to the following 
involutions of first type : $\mu\otimes 1$, $\exp\ i\pi\ad(p_l)$, 
$\exp\ i\pi\ad(\ocirc{p_l})$, $\exp\ i\pi\ad(\ocirc{p_j})$ for $1\leq j\leq 
\frac{l}{2}$ and owing to the equivalence relations of the kind described in 
Example 3, the non-equivalent \aV diagrams for $\frk{g}$ correspond to the 
following involutions of first type : $\rho$, 
$\rho(\mu^2\exp\ i\pi\ad({\ocirc{p_1}})\otimes1)\exp\ i\pi\ad(\ocirc{p_l})$,
$\rho\exp\ i\pi\ad(\ocirc{p_j})$ for $1\leq j\leq \frac{l+1}{2}$. The
non-equivalent \aV diagrams for the twisted affine Kac-Moody Lie 
algebras of type $A_{2l}^{(2)}$, $D_{2r+1}^{(2)}. r\geq 1$, 
$D_{2r}^{(2)}, r\geq 2$, $E_{6}^{(2)}$ and $D_4^{(2)}$ are similarly studied.

It can be easily observed from the Figures 2 and 3 that the count of the 
number of the \aV diagrams corresponding to non-trivial involutions of first 
type matches with the number of almost compact non-compact real forms of 
twisted affine Kac-Mody Lie algebras as given in \cite{M-R}, thereby 
suggesting the existence of a bijective correspondence between the equivalence 
classes of the Vogan diagrams and the isomorphism classes of the almost 
compact real forms of twisted affine Kac-Moody Lie algebras.

\begin{figure}
\begin{tabular}{|l|c|c|} \hline\hline
$\frk{g}$ & 
\multicolumn{1}{c|}{\emph{e-Vogan diagram}} &
\multicolumn{1}{c|}{\emph{involution of first type}}\\ \hline
$A_{2}^{(2)}$ & $\ \underset{\alpha_0}{\stackrel{1}{\bullet}} 
 \,^{=}_{=} \ >\underset{\alpha_{1}}{\stackrel{2}{\circ}}$ & $\mu\otimes1$
\\ \cline{2-3}
& $\quad \underset{\alpha_{0}}{\stackrel{1}{\circ}}  \,^{=}_{=} \ >
  \underset{\alpha_{1}}{\stackrel{2}{\bullet}}$ & $\exp i\pi\ad(\ocirc{p_1})=
\exp i\pi\ad(p_1)$ \\ \hline 
$ A_{2l}^{(2)}$ & $\underset{\alpha_{0}}{\stackrel{1}{\bullet}} 
\Rightarrow \underset{\alpha_{1}}{\stackrel{2}{\circ}}   - \cdots-
  \underset{\alpha_{l-1}}{\stackrel{2}{\circ}} \Rightarrow
  \underset{\alpha_{l}}{\stackrel{2}{\circ}}$ & $\mu\otimes1$\\\cline{2-3}
& $\underset{\alpha_{0}}{\stackrel{1}{\circ}} \Rightarrow
  \underset{\alpha_{1}}{\stackrel{2}{\circ}}   - \cdots-
 \underset{\alpha_{i}}{\stackrel{\overset{\qquad}{2}}{\bullet}}-\cdots
  \underset{\underset{\qquad}{\alpha_{l-1}}}{\stackrel{2}{\circ}} \Rightarrow
  \underset{\alpha_{l}}{\stackrel{2}{\circ}}$ &  
$\begin{array}{c}\exp i\pi\ad(\ocirc{p_i})=
\exp i\pi\ad(p_i),\\ 1\leq i\leq l\end{array}$ \\  \hline
$A_{2l-1}^{(2)}$ $(l \geq 3 )$ &
$\begin{array}{l} \,\, \qquad \stackrel{1}{\bullet} _{\alpha_0}\\
\, \qquad \, \vert \\ 
   \underset{\alpha_{1}}{\stackrel{1}{\circ}} -
     \underset{\alpha{2}}{\circ^{2}} -
   \underset{\alpha_{3}}{\stackrel{2}{\circ}} - \cdots-
   \underset{\alpha_{l-1}}{\stackrel{2}{\circ}} \Leftarrow
   \underset{\alpha_{l}}{\stackrel{1}{\circ}}
 \end{array}$ 
 & $\mu\otimes1$ \\\cline{2-3} 
& $\begin{array}{l} \,\, \qquad \stackrel{1}{\circ} _{\alpha_0}\\
\, \qquad \, \vert \\
   \underset{\alpha_{1}}{\stackrel{1}{\circ}} -
     \underset{\alpha{2}}{\circ^{2}} -
   \underset{\alpha_{3}}{\stackrel{2}{\circ}} - \cdots-
   \underset{\alpha_{l-1}}{\stackrel{2}{\circ}} \Leftarrow
   \underset{\alpha_{l}}{\stackrel{1}{\bullet}}
 \end{array}$ & $(\mu\otimes1)\exp i\pi\ad(\ocirc{p_l}) =\exp i\pi\ad(p_l) $ \\
\cline{2-3}
&$\begin{array}{l} \,\, \qquad \stackrel{1}{\bullet} _{\alpha_0}\\
\, \qquad \, \vert \\
   \underset{\alpha_{1}}{\stackrel{1}{\circ}} -
     \underset{\alpha{2}}{\circ^{2}} -
   \underset{\alpha_{3}}{\stackrel{2}{\circ}} - \cdots-
   \underset{\alpha_{l-1}}{\stackrel{2}{\circ}} \Leftarrow
   \underset{\alpha_{l}}{\stackrel{1}{\bullet}}
 \end{array}$ & $\exp i\pi\ad(\ocirc{p_l})$ \\\cline{2-3}
& $\begin{array}{l} \, \qquad \stackrel{1}{\circ} _{\alpha_0}\\
\, \qquad \, \vert \\
   \underset{\alpha_{1}}{\stackrel{1}{\circ}} -
     \underset{\alpha{2}}{\circ^{2}} -\cdots
   \underset{\alpha_{i}}{\stackrel{2}{\bullet}} - \cdots-
   \underset{\alpha_{l-1}}{\stackrel{2}{\circ}} \Leftarrow
   \underset{\alpha_{l}}{\stackrel{1}{\circ}}
 \end{array}$  & $\exp i\pi\ad(\ocirc{p_i})$, $1\leq i\leq \frac{l}{2}$ \\ 
\hline

$A_{2l-1}^{(2)}$ $(l >2 )$ & 
$\begin{array}{l} \overset{\alpha_0}{_1\circ}   \\
\overset{\diagdown}
{\underset{\diagup} {\updownarrow \quad \circ_{\alpha_{_2}}^2}}
-  \underset{\alpha_{3}}{\stackrel{2}{\circ}} - \cdots-
   \underset{\alpha_{l-1}}{\stackrel{2}{\circ}} \Leftarrow
   \underset{\alpha_{l}}{\stackrel{1}{\circ}} \\
\underset{\alpha_1}{_1\circ} 
\end{array}$ & $\rho$ \\ \cline{2-3}
&$\begin{array}{l} \overset{\alpha_0}{_1\circ}   \\
\overset{\diagdown}
{\underset{\diagup} {\updownarrow \quad \circ_{\alpha_{_2}}^2}}
-  \underset{\alpha_{3}}{\stackrel{2}{\circ}} - \cdots-
   \underset{\alpha_{l-1}}{\stackrel{2}{\circ}} \Leftarrow
   \underset{\alpha_{l}}{\stackrel{1}{\bullet}}\\
\underset{\alpha_1}{_1\circ} 
\end{array}$ & $\rho(\mu^2\exp i\pi\ad(\ocirc{p_1})\otimes1)
\exp i\pi\ad(\ocirc{p_l})$ \\ \cline{2-3}
&  $\begin{array}{l} \overset{\alpha_0}{_1\circ}  \\
\overset{\diagdown}
{\underset{\diagup} {\updownarrow \quad \circ_{\alpha_{_2}}^2}}
-\cdots-  \underset{\alpha_{i}}{\stackrel{2}{\bullet}} - \cdots-
   \underset{\alpha_{l-1}}{\stackrel{2}{\circ}} \Leftarrow
   \underset{\alpha_{l}}{\stackrel{1}{\circ}}\\
\underset{\alpha_1}{_1\circ}\end{array}$ & $\rho\exp i\pi\ad(\ocirc{p_i})$, 
$2\leq i\leq \frac{l+1}{2}$ \\ \hline
\hline\hline \end{tabular}
\caption{Vogan diagrams for affine Kac-Moody Lie algebras of type Aff 2 }
\end{figure}

\begin{figure}
\begin{tabular}{|l|c|c|} 
\hline\hline $\frk{g}$ & \multicolumn{1}{c|}{\emph{e-Vogan diagram}} &
\multicolumn{1}{c|}{\emph{involution of first type}}\\ \hline
$D_{l+1}^{\left(2\right)} \left(l \geq 2 \right)$&
$  \underset{\alpha_{0}}{\stackrel{\overset{\qquad}{1}}{\bullet}}  
\Leftarrow  \underset{\alpha_{1}} {\stackrel{1}{\circ}}  - \cdots- 
\underset{\alpha_{i}} {\stackrel{1}{\circ}}  - \cdots-
  \underset{\alpha_{l-1}}{\stackrel{1}{\circ}}  \Rightarrow
\underset{\underset{\qquad}{\alpha_{l}}}{\stackrel{1}{\circ}}$ &$\mu\otimes1$
\\ \cline{2-3}
&$  \underset{\alpha_{0}}{\stackrel{\overset{\qquad}{1}}{\circ}}  
\Leftarrow  \underset{\alpha_{1}} {\stackrel{1}{\circ}}  - \cdots- 
\underset{\alpha_{i}} {\stackrel{1}{\bullet}}  - \cdots-
  \underset{\alpha_{l-1}}{\stackrel{1}{\circ}}  \Rightarrow
\underset{\underset{\qquad}{\alpha_{l}}}{\stackrel{1}{\circ}}$ &
$\begin{array}{l}\exp i\pi\ad(\ocirc{p_i})\mu\otimes1=\\
\exp i\pi\ad(p_i), 1\leq i\leq \frac{l}{2}\end{array}$ \\ 
\cline{2-3}
&$\underset{\alpha_{0}}{\stackrel{1}{\bullet}}  \Leftarrow
  \underset{\alpha_{1}} {\stackrel{1}{\circ}}  - \cdots- 
\underset{\alpha_{i}} {\stackrel{\overset{\qquad}1}{\bullet}}  - \cdots-
  \underset{\alpha_{l-1}}{\stackrel{1}{\circ}}  \Rightarrow
\underset{\underset{\qquad}{\alpha_{l}}}{\stackrel{1}{\circ}}$ &
$\exp i\pi\ad(\ocirc{p_i})$,$1\leq i\leq l$\\ \hline
$D_{2r+1}^{(2)} \left(r \geq 1 \right)$&
$\begin{array}{l}
\overset{\alpha_{0}}{\circ_1} \Leftarrow
\overset{\quad \alpha_{1}\, \,} {\circ_1} - \cdots- 
\overset{\alpha_{r-1}}{\circ_1}\,\diagdown\\
\updownarrow \, \qquad \, \, \, \updownarrow \qquad 
 \,   \quad  \qquad \updownarrow\quad 
{\alpha_r}\overset{1}{\circ} \quad\\
\underset{\alpha_{2r}}{\circ^1} \Leftarrow 
\underset{\alpha_{2r-1}} {\circ^1} - \cdots- 
\underset{\alpha_{r+1}\,\,}{\circ^1}\diagup
 \end{array}$ & $\rho$ \\ \cline{2-3}
&$\begin{array}{l}
\overset{\alpha_{0}}{\circ_1} \Leftarrow
\overset{\quad \alpha_{1}\, \,} {\circ_1} - \cdots- 
\overset{\alpha_{r-1}}{\circ_1}\,\diagdown\\
\updownarrow \, \qquad \, \, \, \updownarrow \qquad 
 \,   \quad  \qquad \updownarrow\quad 
{\alpha_r}\overset{1}{\bullet} \quad\\
\underset{\alpha_{2r}}{\circ^1} \Leftarrow 
\underset{\alpha_{2r-1}} {\circ^1} - \cdots- 
\underset{\alpha_{r+1}\,\,}{\circ^1}\diagup
\end{array}$ & $\begin{array}{l}\rho\exp i\pi\ad(p_r)=\\
\rho\exp i\pi\ad(\ocirc{p_r})(\mu^2\exp i\pi\ad(\ocirc{p_2r})\otimes1)
\end{array}$ \\ \hline
$D_{2r}^{(2)} \left(r \geq 2 \right)$&
$\begin{array}{c} \overset{\alpha_{0}}{\circ_1} \Leftarrow
\overset{\quad \alpha_{1} \, \, \,} {\circ_1} - \cdots- 
\overset{\alpha_{r-1}}{\circ_1}\\
\updownarrow \, \qquad \,\,\, \updownarrow \qquad 
\,   \quad  \qquad \updownarrow\, \Big)\\
\underset{\alpha_{2r-1}}{\circ^1} \Leftarrow
\underset{\alpha_{2r-2}}{\circ^1} - \cdots- 
\underset{_{\alpha_{r+1}}}{\circ^1}\,  \quad
\end{array}$ & $\rho$ \\ \hline
$E_{6}^{(2)}$ & $ \underset{\alpha_{0}}{\stackrel{\overset{\qquad}{1}}{\circ}}
 -  \underset{\alpha_{1}}{\stackrel{2}{\bullet}}  -
  \underset{\alpha_{2}}{\stackrel{3}{\circ}}  \Leftarrow  
  \underset{\underset{\qquad}{\alpha_{3}}}{\stackrel{2}{\circ}}  -
  \underset{\alpha_{4}}{\stackrel{1}{\circ}} $ & $\begin{array}{l}
\exp i\pi\ad(\ocirc{p_1})\\=\exp i\pi\ad(p_1) \end{array} $\\ \cline{2-3}
& $ \underset{\alpha_{0}}{\stackrel{\overset{\qquad}{1}}{\bullet}}  - 
  \underset{\alpha_{1}}{\stackrel{2}{\circ}}  -
  \underset{\alpha_{2}}{\stackrel{3}{\circ}}  \Leftarrow  
  \underset{\alpha_{3}}{\stackrel{2}{\circ}}  -
  \underset{\underset{\qquad}{\alpha_{4}}}{\stackrel{1}{\bullet}} $& 
$\exp i\pi\ad(\ocirc{p_4}) $
\\ \cline{2-3}
& $ \underset{\alpha_{0}}{\stackrel{\overset{\qquad}{1}}{\circ}}  - 
  \underset{\alpha_{1}}{\stackrel{2}{\circ}}  -
  \underset{\alpha_{2}}{\stackrel{3}{\circ}}  \Leftarrow  
  \underset{\underset{\qquad}{\alpha_{3}}}{\stackrel{2}{\circ}}  -
  \underset{\alpha_{4}}{\stackrel{1}{\bullet}} $& $\begin{array}{l}
(\mu\otimes1)\exp i\pi\ad(\ocirc{p_4})\\=\exp i\pi\ad(p_4)  \end{array}$ \\
\cline{2-3}
& $ \underset{\alpha_{0}}{\stackrel{\overset{\qquad}{1}}{\bullet}}- 
  \underset{\alpha_{1}}{\stackrel{2}{\circ}}  -
  \underset{\alpha_{2}}{\stackrel{3}{\circ}}  \Leftarrow  
  \underset{\underset{\qquad}{\alpha_{3}}}{\stackrel{2}{\circ}}  -
  \underset{\alpha_{4}}{\stackrel{1}{\circ}} $& $\mu\otimes1 $\\ \hline
$D_{4}^{(3)}$ & $\underset{\alpha_{0}}{\stackrel{\overset{\qquad}{1}}{\bullet}}
- \underset{\underset{\qquad}{\alpha_{1}}}{\stackrel{2}{\circ}}   \Lleftarrow
 \underset{\alpha_{2}}{\stackrel{1}{\bullet}} $ & $\exp i\pi\ad(\ocirc{p_2}) $
\\ \cline{2-3}
& $\underset{\alpha_{0}}{\stackrel{\overset{\qquad}{1}}{\bullet}}    -
 \underset{\underset{\qquad}{\alpha_{1}}}{\stackrel{2}{\circ}}   \Lleftarrow
 \underset{\alpha_{2}}{\stackrel{1}{\circ}} $ & $\zeta^3\otimes1 $\\
\hline\hline \end{tabular}
\caption{Vogan diagrams for affine Kac-Moody Lie algebras of type Aff 2 and 3 
(contd.)}
\end{figure}

\newpage
\subsection*{\large{Acknowledgments}}
I thank Dr. Punita Batra for suggesting the problem and providing
me with the references \cite{Bausch} and \cite{R4}.
\bibliographystyle{amsplain}

\begin{bibliography}{10}
\bibliographystyle{amsplain}

\end{bibliography}

\end{document}